\documentclass[12pt]{amsart}
\usepackage{psfig,amssymb,epic,eepic,verbatim}
\input psfig.sty

\newtheorem{theorem}{Theorem}[subsection]

\newtheorem{proposition}[theorem]{Proposition}

\newtheorem{lemma}[theorem]{Lemma}
\newtheorem{lem}[theorem]{Lemma}
\theoremstyle{definition}
\newtheorem{definition}[theorem]{Definition}
\theoremstyle{remark}
\newtheorem{remark}[theorem]{Remark}
\newtheorem{example}[theorem]{Example}
\newcommand{\blem}{\begin{lem} \rm}
\newcommand{\elem}{\end{lem}}
\newcommand{\bml}{\begin{multline}}
\newcommand{\eml}{\end{multline}}
\newcommand\A{\mathcal{A}}
\newcommand\M{\mathcal{M}}

\renewcommand\M{\mathcal{M}}

\renewcommand{\L}{\mathcal{L}}
\renewcommand{\O}{\mathcal{O}}

\newcommand{\N}{\mathbb{N}}
\newcommand{\R}{\mathbb{R}}

\newcommand{\C}{\mathbb{C}}
\newcommand{\cC}{\mathcal{C}}

\newcommand{\Z}{\mathbb{Z}}
\newcommand{\Q}{\mathbb{Q}}
\newcommand{\ddt}{\frac{d}{dt}}

\newcommand{\ddth}{\frac{d}{d \theta}}

\newcommand{\ddu}{\frac{d}{du}}

\newcommand{\ppxi}{\frac{\partial}{\partial \xi}}

\renewcommand{\P}{\mathbb{P}}

\newcommand\lie[1]{\mathfrak{#1}}
\renewcommand{\k}{\lie{k}}
\newcommand{\h}{\lie{h}}
\newcommand{\g}{\lie{g}}

\newcommand{\p}{\lie{p}}

\renewcommand{\t}{\lie{t}}
\newcommand{\Alc}{\lie{A}}

\newcommand{\on}{\operatorname}

\newcommand{\crit}{\on{crit}}
\newcommand{\hol}{\on{hol}}

\newcommand{\s}{\on{s}}

\newcommand{\Aut}{ \on{Aut} }

\newcommand{\Ad}{ \on{Ad} } 
\newcommand{\ad}{ \on{ad} } 
 
\newcommand{\Hol}{ \on{Hol} }

\renewcommand{\ker}{ \on{ker}}
\newcommand{\coker}{ \on{coker}}
\newcommand{\im}{ \on{im}}

\newcommand{\diag}{  \on{diag}}

\newcommand{\ssm}{\kern-.5ex \smallsetminus \kern-.5ex}

\newcommand\dirac{/\kern-1.2ex\partial} 
\newcommand\qu{/\kern-.7ex/} 
\newcommand\lqu{\backslash \kern-.7ex \backslash} 
\newcommand\bs{\backslash}
\newcommand\dr{r_+ \kern-.7ex - \kern-.7ex r_-}
 
\newcommand{\Waff}{W_{\on{aff}}} 

\newcommand{\labell}\label

\renewcommand{\d}{{\mbox{d}}}
\newcommand{\ol}{\overline}

\newcommand\Phinv{\Phi^{-1}}

\newcommand\eps{\epsilon}
\newcommand\Om{\Omega}
\newcommand\om{\omega}

\newcommand{\f}{\frac}

\newcommand{\hh}{{\f{1}{2}}}
\newcommand{\qq}{{\f{1}{4}}}
\newcommand{\thh}{{\f{3}{2}}}
\newcommand{\ti}{\tilde}

\renewcommand{\ss}{\on{ss}}
\newcommand\Tr{\on{Tr}}
\newcommand\mE{\mathcal{E}}

\newcommand\Map{\on{Map}}

\newcommand\ul{\underline}

\renewcommand\Im{\on{Im}}
\newcommand\Ker{\on{Ker}}
\newcommand\grad{\on{grad}}

\newcommand\bdefn{\begin{definition}}
\newcommand\edefn{\end{definition}}
\newcommand\bea{\begin{eqnarray*}}
\newcommand\eea{\end{eqnarray*}}
\newcommand\bcv{\left[ \begin{array}{r} }
\newcommand\ecv{\end{array} \right] }

\newcommand\bma{\left[ \begin{array} }
\newcommand\ema{\end{array} \right]}
\newcommand\ben{\begin{enumerate}}
\newcommand\een{\end{enumerate}}
\newcommand\bex{\begin{example}}
\newcommand\eex{\end{example}}
\newcommand\sx{*\kern-.5ex_X}

\def\mathunderaccent#1{\let\theaccent#1\mathpalette\putaccentunder}
\def\putaccentunder#1#2{\oalign{$#1#2$\crcr\hidewidth \vbox
to.2ex{\hbox{$#1\theaccent{}$}\vss}\hidewidth}}

\begin{document}

\title[The heat flow on the moduli space of framed bundles]{The Yang-Mills heat flow on the moduli space of framed bundles on
  a surface}

\author{Christopher T. Woodward} 
\thanks{Partially supported by NSF
grants DMS9971357 and DMS0093647 } 
\thanks{An earlier version was
titled {\em A Kirwan-Ness stratification for loop groups}} 
\address{Mathematics-Hill Center,
Rutgers University, 110 Frelinghuysen Road, Piscataway NJ 08854-8019,
U.S.A.}  \email{ctw@math.rutgers.edu}

\begin{abstract} 
We study the analog of the Yang-Mills heat flow on the moduli space of
framed bundles on a cut surface.  Existence and convergence of the
heat flow give a stratification of Morse type invariant under the
action of the loop group. We use the stratification to prove versions
of K\"ahler quantization commutes with reduction and Kirwan
surjectivity.
\end{abstract}

\maketitle

\section{Introduction}

Let $K$ be a compact, $1$-connected Lie group with complexification
$G$ and Lie algebra $\k$, and $\ol{X}$ a compact, connected Riemann
surface.  The moduli space $\M(\ol{X})$ of isomorphism classes of flat
$K$-bundles on $\ol{X}$ is homeomorphic to the moduli space of
grade-equivalence classes of semistable $G$-bundles, by the theorems
of Narasimhan-Seshadri \cite{ns:ho} and Ramanathan
\cite{ra:th,ra:th2}.

$\M(\ol{X})$ has two presentations as an infinite dimensional quotient
which can be used to study its cohomology.  The first presentation was
introduced by Atiyah and Bott \cite{at:mo} and is rather well
understood.  Let $\A(\ol{X})$ denote the affine space of connections
on the trivial $K$-bundle over $X$, with symplectic structure induced
by a choice of metric on $\k$.  The group $K(\ol{X})$ of gauge
transformations acts symplectically on $\A(\ol{X})$ with moment map
given by the curvature, and the symplectic quotient is $\M(\ol{X})$.
In the holomorphic description
$$ \M(\ol{X}) \cong G(\ol{X}) \lqu \A(\ol{X}) $$
where the symbol $\lqu$ means the quotient of the semistable locus.
Atiyah and Bott used the stratification of $\A(\ol{X})$ into
Harder-Narasimhan types to compute the Betti numbers of $\M(\ol{X})$;
they conjectured that the stratification is identical to the
stratification into stable manifolds for the gradient flow of minus
the Yang-Mills functional.  This was proved by Donaldson \cite{do:ne}
and Daskalopolous \cite{da:st}.  R{\aa}de \cite{ra:ym} proved that the
gradient flow converges.

The second presentation has older origins (Weil's double coset
construction) but the related analysis has been less studied.  Let $S
\subset \ol{X}$ be an embedded circle, and $X$ the Riemann surface
with boundary obtained by cutting $\ol{X}$ along $S$.  The Yang-Mills
heat flow on the space $\A(X)$ was studied by Donaldson \cite{do:bv},
who obtained an analog of the Narasimhan-Seshadri theorem: The moduli
space $\M(X)$ of flat $K$-bundles with framings on the boundary is
diffeomorphic to $G(\partial X)/G_{\hol}(X)$, where $G(\partial X)=
\Map(\partial X,G)$ and $G_{\hol}(X)$ denotes the subgroup of
$G(\partial X)$ consisting of loops that extend holomorphically over
the interior.  The loop group $K(S)$ acts symplectically on $\M(X)$
with moment map given by the difference of the restriction to the two
boundary components, and $\M(\ol{X})$ is homeomorphic to the
symplectic quotient.  In the holomorphic description
$$ \M(\ol{X}) \cong G(S) \lqu ( G(\partial
X)/G_{\hol}(X)) .$$
In recent years, this presentation has become more popular because of
its connection with conformal field theory and the Verlinde formulas
\cite{be:cb}, \cite{fa:vl}, \cite{ku:in}, \cite{te:bo}.  Here the
circle $S$ is assumed to bound a disk, so that $X$ (in its algebraic
manifestation) becomes a punctured curve union a formal disk.
Surfaces with boundary do not fit into the algebraic framework.

In this paper we study the analog of the Yang-Mills heat flow in the
second presentation, namely the gradient flow of minus the square of
the moment map for the loop group for an arbitrary embedded circle $S$
in $\ol{X}$.  We show that the analog of R{\aa}de's result holds: The
gradient flow exists for all times and converges to a critical point.
Although the evolution equation itself is not pseudo-differential, its
restriction to the boundary is a (non-linear) heat equation involving
the Dirichlet-to-Neumann operator associated to the connection.
Calder\'on observed that this is an elliptic pseudodifferential
operator.  Because pull-back to the boundary is Fredholm on the space
of harmonic forms, we are ``up to finite dimensions'' in the same
situation as for the first presentation, except that the moduli space
of framed bundles is not affine.

This analysis implies that $\M(X)$ admits a stratification into stable
manifolds for minus the gradient flow.  By definition, the stable
manifold for the zero locus of the moment map is the semistable locus.
The other strata are complex submanifolds of finite codimension, and
the number of strata of each codimension is finite.  Using the
stratification, we obtain several cohomological applications which
extend known results beyond the case that $S$ bounds a disk.  The
first, which was the motivation for the paper, is a K\"ahler
``quantization commutes with reduction'' theorem, similar to that of
Guillemin-Sternberg \cite{gu:ge} in the finite dimensional case.  This
is an instance of Segal's composition axiom for the Wess-Zumino-Witten
conformal field theory, see the recent thesis of H. Posthuma
\cite{po:th}.  In the case $S$ bounds a disk, the algebraic version is
due to Beauville-Laszlo \cite{be:cb}, Kumar-Narasimhan-Ramanathan
\cite{ku:in}, and Laszlo-Sorger \cite{la:li}.  We had in mind the case
that $X$ is a four-pointed projective line and $S$ separates the four
points into two groups of two, which is related to asymptotics of $6j$
symbols \cite{ro:qu}, \cite{ta:6j}.  The second application is a a
surjectivity result for the equivariant cohomology with rational
coefficients, similar to that of Kirwan \cite{ki:co}.  In the case
that $S$ bounds a disk this was recently proved by Bott, Tolman, and
Weitsman \cite{btw:ta}.

Two appendices contain a review of Sobolev spaces, and proof of
convergence of the gradient flow of minus the norm-square of the
moment map in finite dimensions.

\section{Background on connections on a circle}

The following is contained in Pressley-Segal \cite{ps:lg}, except for
smooth maps.  Let $S$ be a circle, that is, a connected one-manifold.
For any $s > 0$, the group
$$K(S)_{s+ \hh} := \Map(S,K)_{s + \hh}$$ 
of free loops of Sobolev class $s + \hh$ acts on the space $\A(S)_{s -
\hh}$ of connections on the trivial bundle $S \times K$.  Any
connection differs from the trivial connection by a $\k$-valued
one-form; using the trivial connection as a base point we identify
$$\A(S)_{s - \hh} \to \Omega^1(S;\k)_{s - \hh} .$$  
For any $ s > r > 0$ inclusion defines a bijection
$$ K(S)_{r - \hh} \backslash \A(S)_{r + \hh} \to K(S)_{s - \hh}
\backslash \A(S)_{s + \hh} .$$
For $s > 2$, there is a smooth holonomy map
$$\on{\Hol}: \ \A(S)_{s - \hh} \to K$$
depending on the choice of base point $*$ in $S$; the assumption $s >
2$ implies that $A$ is $C^1$ which guarantees existence of a solution
to the parallel transport equation.  For $s > 0$ and $A \in \A(S)_{s -
\hh}$ the stabilizer $K(S)_A$ is a compact, connected Lie group.  For
$s > 2$, $K(S)_{s+\hh,A}$ is isomorphic to the centralizer of the
holonomy $\on{Hol}(A)$ via the map
$$K(S)_{s+\hh} \to K, \ \ k \mapsto k(*) .$$
For $s > 0$ there are bijections
$$K_*(S)_{s + \hh} \backslash \A(S)_{s -\hh} \to K, \ \ \ K(S)_{s +
\hh} \backslash \A(S)_{ s- \hh} \to \Ad(K) \backslash K,$$
which for $s > 2$ are given by taking the holonomy, resp. conjugacy
class of the holonomy of the connection.

The orbits of $K(S)_{s + \hh}$ on $\A(S)_{s - \hh}$ can be
parametrized by the Weyl alcove as follows.  Let $\Lambda$ denote the
coweight lattice of $T$ and
$$ \Waff := W \rtimes \Lambda $$
the affine Weyl group.  The action of $\Waff$ on the Cartan subalgebra
$\t$ has fundamental domain
$$ \Alc := \{ \xi \in \t_+, \ \alpha_0(\xi) \leq 1 \} $$
where $\t_+$ denotes the positive chamber and $\alpha_0$ the highest
root.  Inclusion and exponentiation define bijections
$$
 \Alc \to \Waff \backslash \t \to W \backslash T \to \Ad(K)
\backslash K  $$
and so for $s > 0$ we have a bijection
\begin{equation} \label{circle}
 \Alc \to K(S)_{s + \hh} \backslash \A(S)_{s- \hh}. \end{equation}

\section{Background on connections on a surface}

\label{surface}

Let $X$ be a compact, connected, oriented surface.  Since $K$ is
simply-connected, any principal $K$-bundle is isomorphic to the
trivial $X \times K$.  Let $A$ be a connection on $X \times K$.  For
any $K$-representation $V$, let
$$ \d_A(V): \ \Omega^0(X;V)_{s+1} \to \Omega^1(X;V)_{s} \to
\Omega^2(X;V)_{s-1} $$
denote the associated covariant derivative.  Let $\d_A := \d_A(\k)$
denote the covariant derivative for the adjoint representation.  $A$
is flat if and only if $\d_A^2 = 0$.  Let $\A(X)_{s}$ denote the
affine space of connections on $X \times K$ of Sobolev class $s$.
Using the trivial connection as base point we may identify
$$\A(X)_s \to \Omega^1(X;\k)_s .$$
For $s > 0$ the gauge group 
$$ K(X)_{s+1} := \Map(X,K)_{s+1}$$ 
is a Banach Lie group and acts on $\A(X)_s$ by the formula
$$ k \cdot A = \Ad(k)A + k d(k^{-1}) = \Ad(k) A - \d k \, k^{-1} $$
in any faithful matrix representation of $K$.  It has Lie algebra
$$ \k(X)_{s+1} := \Omega^0(X;\k)_{s + 1} .$$
The generating vector fields for the action of $K(X)_{s+1}$ on
$\A(X)_s$ are
$$ \xi_{\A(X)} (A) = - \d_A \xi, \ \ \xi \in \k(X)_{s+1} .$$
In particular, the Lie algebra $\k(X)_A$ of the stabilizer $K(X)_A$ of
$A$ is 
$$\k(X)_A = \ker (\d_A | \Omega^0(X;\k)).$$

Suppose $X$ is equipped with a complex structure.  The map
$$d + \ad(A) \mapsto \ol{\partial}_{\alpha} := \ol{\partial} +
\ad(\alpha)$$ 
defines a one-to-one correspondence between covariant derivatives and
holomorphic covariant derivatives
$$\ol{\partial}_{\alpha}: \ \Omega^0(X;\g) \mapsto
\Omega^{0,1}(X;\g)$$
satisfying the holomorphic Leibniz rule $\ol{\partial}_{\alpha}(fs) =
(\ol{\partial} f)s + f \ol{\partial}_{\alpha}s$.  $G(X)$ acts on the
space of holomorphic covariant derivatives by conjugation, and
therefore on the space of $\g$-valued $(0,1)$-forms by
$$ g \cdot \alpha = \Ad(g) \alpha - (\ol{\partial}g) g^{-1} .$$ 
This formula extends to a holomorphic action of $G(X)_{s+1}$ on
$\A(X)_s$.

Using an invariant inner product on $\k$ we define a weakly symplectic
form (that is, a closed $2$-form that defines an injection $T\A(X)_s
\to T^*\A(X)_s$) on $\A(X)_s$ for $s \ge 0$
$$ \omega_{\A(X)}: \ (a_1,a_2) \mapsto \int_X \Tr(a_1 \wedge a_2) .$$
In the case that the boundary of $X$ is empty, the action of
$K(X)_{s+1}$ is Hamiltonian with moment map given by the curvature
\cite{at:mo}
$$ \M(X)_s \to \Omega^2(X;\k)_{s+1}, \ \ A \mapsto F_A .$$
Let $\A_\flat(X)_s$ denote the subspace of flat connections, 
$$ \A_\flat(X)_s := \{ A \in \A(X)_s, \ \ F_A = 0 \} .$$  
The symplectic quotient
$$ \M(X)_s = K(X)_{s+1} \lqu \A(X)_s := K(X)_{s+1} \backslash
\A_\flat(X)_s $$
is the {\em moduli space of flat bundles} on $X$.

In the case $X$ has boundary, the moment map picks up an additional
term \cite{at:kn}, \cite{do:bv}, \cite{me:lo}
\begin{equation} \label{moment}
\A(X)_s \to \Omega^2(X;\k)_{s-1} \oplus \Omega^1(\partial X;\k)_{s-1},
\ \ \ A \mapsto (F_A, r_{\partial X} A) .\end{equation}
That is, for all $\xi \in \k(X)_{s+1}$
$$ \iota(\xi_{\A(X)}) \omega_{\A(X)} = - d \left( \int_X \Tr(F_A
\wedge \xi) - \int_{\partial X} \Tr(r_{\partial X} A \wedge \xi)
\right) .$$
Let $K_\partial(X)_{s+1}$ the subgroup fixing a framing on the
boundary,
$$ K_\partial(X)_{s+1} = \{ k \in K(X)_{s+1}, \ k |_{\partial X} = 1
\} .$$
For $s> 0$ there is an exact sequence of Banach Lie groups
$$1 \to K_\partial (X)_{s+1} \to K(X)_{s+1} \to K(\partial X)_{s+\hh}
\to 1.
$$
Surjectivity of the third map follows from triviality of $\pi_1(K)$
and the properties of the extension operator \ref{sobolev}
\eqref{extension}.  The moment map for $K_\partial(X)_{s+1}$ is the
curvature and the symplectic quotient
$$ \M(X)_s = K_\partial(X)_{s+1} \lqu \A(X)_s  :=  
K_\partial(X)_{s+1} \bs \A_{\flat}(X)_s 
$$
is the {\em moduli space of framed flat bundles} on $X$.

Charts for $\M(X)_s$ are constructed from local slices for the gauge
action as follows.  First note that $d_A$ is surjective: $\coker (d_A:
| \Omega^1(X;\k)_s)$ is isomorphic to $ \ker (d_A |
\Omega^0(X,\partial X;\k)_{1-s})$.  By elliptic regularity this is
contained in $ \ker (d_A | \Omega^0(X,\partial X;\k)_1)$ and therefore
trivial.  Hence there exists a right inverse $\d_A^{-1}:
\Omega^2(X;\k)_{s-2} \to \Omega^1(X;\k)_{s-1}$ for $\d_A$, depending
continuously on $A$; for $A$ smooth this inverse is
pseudodifferential.  By the implicit function theorem, there exists a
constant $\eps$ depending only on $\Vert \d_A^{-1} \Vert$, open
neighborhoods of $A$, resp. $0$
\begin{equation} \label{UA}
 U_A \subset \{ A + a \in \A(X)_s, \ \ F_{A + a} = 0, \ \ \ \d_A^* a
= 0 \} \end{equation}
\begin{equation} \label{VA}
 V_A \subset \{ a \in \Omega^1(X;\k)_{s}, \ \ \d_A a = 0, \ \ \
\d_A^* a = 0 \} \end{equation}
such that $V_A $ is an $\eps$-ball around $0$, and a smooth map 
$$S: \ V_A \to \Omega^2(X;\k)_{s-1}$$ 
such that
\begin{equation} \label{Feq}
 F_{A + (I + \d_A^{-1}S) a} = Sa + \hh [a + \d_A^{-1}Sa,a +
 \d_A^{-1}Sa] =0.
\end{equation}
Define 
$$ \varphi_A: \ V_A \to U_A, \ \ a \mapsto A + (I +
\d_A^{-1}S)a.$$

The following lemma summarizes the basic properties
of $\M(X)_s$:

\blem \label{surf} \ben
\item For any $ s> 0$, $\M(X)_s$ is a smooth Banach manifold.
\item The size of the slice $U_A$ at $A \in \A(X)_s$ depends only on
the operator norms $\Vert \d_A \Vert$ and $\Vert \d_A^{-1} \Vert$.
\item For any $r > s > 0$, the inclusion $\A_\flat(X)_r \to
\A_\flat(X)_s$ induces a bijection
$$K(\partial
X)_{r + 1} \backslash (K(\partial X)_{s + 1} \times \M(X)_r) \to
\M(X)_s .$$
\item If $s > \qq$ restriction to the boundary
$$ \M(X)_s \to \Omega^1(\partial X, \k)_{s - \hh}, \ \ [A] \mapsto
r_{\partial X} A $$
is a proper moment map for the action of $K(\partial X)_{s + \hh}$.
\item 
 \label{holonomy} For any $ s>2$, $\M(X)_s$ is diffeomorphic to a
fiber product 
$$K^{2(g + b - 1)} \times_{K^b} \Omega^1(\partial
X;\k)_{s - \hh}  .$$  
The maps are described in the proof below.
\item For any $s > \qq$, the quotient $K(\partial X)_{s +1} \backslash
\M(X)_s$ is compact.

\item For any $s > 0$, the Hodge star 
$$\sx: \ker \d_A \oplus \d_A^* \to \ker \d_A \oplus \d_A^* $$ 
defines an almost complex structure on $\M(X)_s$.

\item For any $s > 0$, $\M(X)_s $ is diffeomorphic to $G(\partial
X)_{s+\hh}/G_{\hol}(X)_{s+\hh}$.

\item For any $s > 0$, the almost complex structure $\sx$ is
integrable, that is, there exist charts for which the transition maps
are holomorphic.  \een \elem

\begin{proof}
The proofs of (a) and (c) are somewhat standard, as in \cite{do:fo},
and left to the reader.  (b) follows since the size of the ball in
$V_A$ depends only on $\Vert \d_A^{-1} \Vert, \Vert \d_A \Vert$.  (d)
If $s > \hh$ holds then the trace is a continuous linear map $
r_{\partial X} : \ \Om^1(X;\k)_s \to \Omega^1(\partial X;\k)_{s- \hh};
$.  It follows from \eqref{moment} that $r_{\partial X}$ is a moment
map.  Therefore, the problem is to establish the lemma in the case
$\qq < s < \hh$.  To get improved regularity we wish to show that for
$A \in \A(X)_s$ the trace map $ r_{\partial X} \circ \varphi_A : \
V_{A,s}^\eps \to \Omega^1(\partial X;\k)_{s- \hh}$ is smooth.  Since
points $[A]$ with $A$ smooth are dense, we may assume that $A$ is
smooth.  Consider the terms in \eqref{Feq}.  The element $a$ satisfies
$(\d_A \oplus \d_A^*)a = 0$, hence $(d \oplus d^*)a = - (\Ad(A)+
\Ad(A)^*)a$. The right hand side has class $2s -1 > - \hh$, which
implies that the trace of $a$ is well-defined \cite[Theorem
13.8]{bo:el}.  The non-linear term $\hh [a + \d_A^{-1}Sa,a +
\d_A^{-1}Sa] $, and therefore also $Sa$, has class $2s - 1 > - \hh.$
Therefore, $\d_A^{-1} Sa$ is class $2s > \hh$ which implies that the
its trace is also defined.  (e) Let $*_1,\ldots,*_b$ be base points on
the boundary components, and let $d_1,\ldots,d_b$ be the paths around
the boundary.  Choose paths $a_1,b_1,\ldots,a_g,b_g$ from $*_1$ to
$*_1$ and $c_1,\ldots,c_{b-1}$ from $*_1$ to $*_2,\ldots,*_{b-1}$ so
that the fundamental group of $X$ is freely generated by $a_i,b_i$, $i
= 1, \ldots,g$ and $\Ad(c_j) d_j, j = 1,\ldots,b-1$.  Suppose $s> 2$.
Then the holonomies $A_i \in K$ around $a_i$, etc.  are well-defined,
and there are smooth maps
$$ \on{Hol}: \ \M(X)_s \to K^{2(g + b - 1)}, \ \ \ \ \ \ \
r_{\partial X}:\ \M(X)_s \to \Omega^1(\partial X;\k)_{s - \hh} .$$
given by $(a,b,c,d)$ and the boundary values.  The remainder of the
proof is the same as in \cite[Theorem 3.2]{me:lo}.  For $s> 2$,
properness of $r_{\partial X}$ follows from compactness of $K^{2g}$ in
the holonomy description.  The extension to $s > \qq$ follows from the
symplectic cross-section theorem: For any face $\sigma^b$ of $\Alc^b$,
let $\Alc^b_\sigma$ denote the open subset of $\Alc^b$ obtained by
removing all faces $\tau$ whose closure $\ol{\tau}$ does not contain
$\sigma$.  Let $K(\partial X)_\sigma$ denote the stabilizer of any
point in $\sigma$.  This is a compact, connected subgroup of
$K(\partial X)_s$, independent of the choice of $s$.  Then
\cite[Section 4.2]{me:lo} $ \Phinv(K(\partial X)_\sigma \Alc_\sigma)$
is a finite dimensional symplectic submanifold of $\M(X)_s$ and there
is a diffeomorphism
$$ \Phinv(K(\partial X)_{s+\hh} \Alc_\sigma) \to K(\partial X)_{s+\hh} 
\times_{K(\partial X)_\sigma} \Phinv(K(\partial X)_\sigma \Alc_\sigma)
.$$
Any compact subset can be written as a finite union of subsets of
$K(\partial X)_s \Alc_\sigma$, as $\sigma$ ranges over faces of
$\Alc^b$.  Therefore it suffices to show that the map
$$    \Phinv(K(\partial X)_s \Alc_\sigma) \to K(\partial X)_s
\Alc_\sigma 
$$ 
is proper.  But this follows from properness of
$$ K(\partial X)_s \times \Phinv(K(\partial X)_\sigma \Alc_\sigma) \to
K(\partial X)_s \times K(\partial X)_\sigma \Alc_\sigma $$
and the fact that $K(\partial X)_\sigma$ is compact.  (f) By (c),
compactness for $r >2$ implies compactness for $s > \qq$.  (g) In
general $*_X^2 = (-1)^{d(\dim(X)-d)}$ on $\Omega^d(X,\k)$.  In this
case $d = 1$ so $*_X^2 = -1$.  (h) For any $\alpha \in
\Omega^{0,1}(X;\g)$, the operator $\ol{\partial}_\alpha $ is
surjective, see e.g. \cite[Chapter 9]{bo:el}.  By the implicit
function theorem, a neighborhood of the identity in $G(X)_{s+1}$ maps
diffeomorphically onto a neighborhood of $\alpha$ in $\A(X)_s$.  It
follows that $G(X)_{s+1}$ acts transitively on $\A(X)_s$, and so
$G(\partial X)_{s+\hh}$ acts transitively on $\M(X)_s$.  The
stabilizer of the trivial connection is $G_{\hol}(X)_{s + \hh}$; it
follows that there is a homeomorphism (in fact, a diffeomorphism of
Banach manifolds) $ \M(X)_s \to G(\partial X)_{s+ \hh}/G_{\hol}(X)_{s
+ \hh}. $ (i) follows from the description in (i), since
$G_{\hol}(X)_{s+\hh}$ is a complex Banach subgroup of $G(\partial
X)_{s + \hh}$. \end{proof}

A {\em marking} is an element $\mu \in \Alc$.  If $\mu_1,\ldots,\mu_b
\in \Alc$ then we define the {\em moduli space of flat bundles with
fixed holonomies}
$$ \M(X;\mu_1,\ldots,\mu_b) = K(\partial X)_{s + \hh} \backslash
r_{\partial X}^{-1}( \O_1 \times \ldots \times \O_b) $$
where $\O_1,\ldots, \O_b$ are the orbits corresponding to
$\mu_1,\ldots,\mu_b$ in \eqref{circle}.  

\subsection{The determinant/Chern-Simons line bundle}

This is a Hermitian line bundle with connection
$$(\L(X)_s,\nabla) \to (\M(X)_s,\omega)$$ 
whose curvature is $- 2 \pi i \omega$, cf. \cite[Section 2]{wi:vl}.
It may be constructed by symplectic reduction as follows \cite{ra:cs},
\cite[Section 3.3]{me:lo}.  The trivial line bundle $\A(X)_s \times
\C$ with connection 1-form
$$ T_A \A(X)_s \to \R, \ \ a \mapsto \hh \int_X \Tr( a \wedge A) .$$
has curvature equal to $- 2\pi i  \om_{\A}$.  The central $U(1)$-extension
$\widehat{K(X)}_{s+1}$ defined by the cocycle
$$ (k_1,k_2) \mapsto \exp \left( \pi i \int_X \Tr( k_1^{-1} \, \d k_1
\wedge \d k_2 \,k_2^{-1}) \right) $$
($k_1^{-1} \, \d k_1$, resp. $\d k_2 \, k_2^{-1}$ are the pull-backs
of the left, resp. right Maurer-Cartan forms on $K$) acts on $\A(X)_s
\times \C$ by connection preserving automorphisms by the formula
$$ (k,z) \cdot (A,w) = (k \cdot A, \, \exp \left( \pi i \int_X 
\Tr(k^{-1} \, \d k \wedge A) \right ) \, z \, w ) .$$
On the Lie algebra level, $\widehat{\k(X)}_{s+1}$ is the central
$\R$-extension of $\k(X)_{s+1}$ defined by the cocycle
\begin{equation} \label{liecocycle}
 (\xi_1,\xi_2) \mapsto \int_X \Tr(d \xi_1 \wedge d \xi_2) = 
\int_{\partial X} \Tr(\xi_1 \, d \xi_2) .\end{equation}
One may use the Chern-Simons three-form to trivialize the restriction
$\widehat{K(X)}_{s+1}$ to $K_\partial(X)_{s+1}$ \cite[Section
3.3]{me:lo}.  The quotient
$$ \widehat{K(\partial X)}_{s + \hh} :=
\widehat{K(X)}_{s+1}/K_\partial(X)_{s+1} $$
is the unique central $U(1)$-extension of $K(\partial X)_{s + \hh}$
defined by the Lie algebra cocycle \eqref{liecocycle}.  Define the
pre-quantum line bundle $\L(X)_s$ by
$$ \L(X)_s = K_\partial(X)_{s+1} \lqu (\A(X)_s \times \C) :=
K_\partial(X)_{s+1} \backslash (A_\flat(X)_s \times \C) .$$
The products $U_A \times \C$ are local slices for the
$K_\partial(X)_{s+1}$- action, and equip $\L(X)_s$ with the structure
of a $\widehat{K(\partial X)}_{s + \hh} $-equivariant Hermitian line
bundle with connection.

The total space $\L(X)_s$ has an almost complex structure $J_\L$
determined by the connection $\nabla$ and the almost complex structure
on $\M(X)_s$, derived from the splitting 
$$T\L(X)_s \cong \pi^* T\M(X)_s \oplus \ul{\C} ,$$
where $\ul{\C}$ denotes the trivial line bundle.

\begin{lemma}  For any $s > 0$, the almost complex structure $J_\L$ is integrable, that is, there exist
local trivializations for which the transition maps for $\L$ are
holomorphic.
\end{lemma}

\begin{proof} Define $\widehat{G(\partial X)}_{s + \hh}$
to be the pull-back of the central $\C^*$-extension
$$ \C^* \to \Aut(\L(X)_s,J_\L) \to \Aut(\M(X)_s,J_\M) $$
under the map $G(\partial X)_{s + \hh} \to \Aut(\M(X)_s,J_\M)$.
Uniqueness of the central extension with a given cocycle implies that
$\widehat{G(\partial X)}_{s + \hh}$ is the {\em basic} central
$\C^*$-extension of ${G(\partial X)}_{s + \hh}$; by \cite[Chapter
6]{ps:lg}, $\widehat{G(\partial X)}_{s + \hh}$ is a complex Banach Lie
group.  Since the action of $G(\partial X)_{s + \hh}$ on $\M(X)_s$ is
transitive, so is the action of $\widehat{G(\partial X)_{s+\hh}}$ on
$\L(X)_s$.  Fix as base point the trivial connection $[0] \in
\M(X)_s$, and let $\L(X)_{[0],s}$ denote the fiber.  The map
$$(\widehat{G(\partial X)_{s+\hh}} \times
\L(X)_{[0],s})/\widehat{G_{\hol}(X)}_{s+\hh} \to \L(X)_s, \ \
[\hat{g},z] \mapsto \hat{g} z $$
is a diffeomorphism preserving the almost complex structure.  Since
$\widehat{G(\partial X)_{s+\hh}}$, $\widehat{G_{\hol}(X)}_{s+\hh} $ are
complex Banach Lie groups, the almost complex structure on the total
space of $\L(X)_s$ is integrable.  Local holomorphic triviality
follows from the existence of local slices.
\end{proof}

\subsection{Gluing equals reduction}

Let $\ol{X}$ be a compact, connected Riemann surface, $S \subset
\ol{X}$ an embedded circle, and $X$ the Riemann surface obtained by
cutting $\ol{X}$ along $S$.  

Let $\pi: X \to \ol{X}$ denote the gluing map, and $S_\pm$ be the
component of $\pi^{-1}(S)$ whose orientation agrees (resp. is the
opposite) of the orientation on $S$.  Let $\pi_\pm$ denote the
restriction of $\pi$ to $S_\pm$.  Consider the diagonal embedding
$$ \delta: \ K(S)_{s + \hh} \to K(\partial X)_{s + \hh}, \ \ \ \ k
\mapsto (\pi_+^* k, \pi_-^* k) .$$
We denote by $r_\pm$ the pull-back (restriction to the boundary,
a.k.a. trace map)
$$ r_\pm: \ \Omega^*(X,\k)_s \to \Omega^*(S_\pm,\k)_{s - \hh} .$$
The moment map for $K(S)_{s + \hh}$ is
$$ \Phi: \ \M(X)_s \to \Omega^1(S;\k)_{s - \hh}, \ \ \ [A] \mapsto
(r_- - r_+) A. $$
The reason for the minus sign is that the identification $S \to S_-$
is orientation reversing.  The map $[A] \mapsto [\pi_X^*A]$ induces a
homeomorphism
$$\M(\ol{X}) \to K(S)_{s+1} \lqu \M(X)_s  ;$$
in fact, an isomorphism of K\"ahler symplectic orbifolds on the
quotient of the subset of $\Phinv(0)$ on which $K(S)_{s+1}$ with only
finite stabilizers \cite[Theorem 3.5]{me:lo}.  

\subsection{Estimates on the sizes of slices}

In this section, we bound the size of gauge slices (charts for
$\M(X)$) from below.  The notation $ < c$ means less than a universal
constant, $< c(R)$ means less than a constant depending on $R$, where
$\hh \Vert (r_+ - r_-)A \Vert_0^2 < R$.

\begin{proposition} 
\label{chart} For any $s \ge \hh$ and  $A \in \A_{\flat,s}$ there
exists a gauge transformation $k \in K(X)_{s+1}$ with $r_+ k = r_- k$
such the open ball $V_{k \cdot A}$ defined in \eqref{VA} has radius is
bounded from below by $c(R)$.\end{proposition}

Since $K(\partial X)_{s+\hh} \backslash \M(X)_s$ is compact, there
exists a compact subset $\A_{\circ}(X)$ of $ \A(X)_s$ such that any
element of $\M(X)_s$ may be represented by an element of
$\A_{\circ}(X)$ up to gauge transform, i.e.
$$ \M(X)_s = \{ k \cdot [A_\circ], \ k \in K(\partial X)_{s + 1}, \ \
A_\circ \in \A_{\circ}(X)_s \}.$$
Let $[A] \in \M(X)_s$ and let $k_\circ \in K(X)_{s+1}$ and $A_\circ \in
\A_{\circ}(X)$ be such that $[k_\circ \cdot A_\circ] = [A]$.  Consider the
operators
$$ \d_A: \Omega^0(X;\k)_{s+1} \to
\Omega^1(X;\k)_s, \ \ \ \d_A^{-1}: \Omega^2(X;\k)_{s-1} \to
\Omega^1(X;\k)_s .$$ 
\blem \label{normlemma} The norms of
$d_A, d_A^{-1}$ can be bounded by a constant depending only on the
norms of $k_\circ$, $k_\circ^{-1}$. \elem

\begin{proof}  We have 
$$
\d_A^{-1} = d_{k_\circ \cdot A_\circ}^{-1} = \Ad(k_\circ ) \circ
d^{-1}_{A_0} \circ \Ad(k_\circ^{-1})  .$$
Since $\A_{\circ}$ is compact, $\Vert \d_A^{-1} \Vert , \Vert \d_A
\Vert $ are bounded on $\A_{\circ}$.  The claim follows. \end{proof}

\blem
\label{pathlemma} 
Let $s \ge \hh$ and let $A \in \A(X)_s$ be a flat connection. There is
a differentiable path $k_t \in K(X)_{s+1}$ such that $A_t = k_t A_0$
satisfies 
\begin{enumerate}
\item  $A_0 \in \A(X)_\circ$, 
\item  $A_1 = k \cdot A$ for some $k \in K(X)_{s+1}$ with $r_{\partial X} k \in K(S)_{s + \hh}$ and
\item  for all $t \in [0,1]$ 
$$\Vert k_t \Vert_{\thh} < c(R), \ \ \ \Vert k_t^{-1} \Vert_{\thh} <
c(R) , \ \ \  \Vert A_t \Vert_{\hh} < c(R), \ \
\text{and} \  \ \Vert \ddt A_t \Vert_{\hh} <
c(R) .$$
\end{enumerate}
\elem

\begin{proof}  Choose $k \in K(X)_2$ with 
$r_+ k = r_- k$ so that $r_+ k \cdot A \in *_S \Alc .$ Let $c =
\sup_{\xi \in \Alc} \Vert *_S \xi \Vert_0 $.  Then
$$ \Vert r_+ k \cdot A \Vert_0 < c \ \ \text{and} \  \ \Vert r_- k \cdot A \Vert_0
\leq \Vert r_+ k \cdot A \Vert_0 + \Vert \Phi(k \cdot A) \Vert_0 <
c(R) .$$
Replacing $A$ with $k \cdot A$, we may assume $\Vert r_{\partial X} A
\Vert_0 < c(R)$.  Let $A_\circ \in \A_\circ(X)_s$ be a flat connection
gauge equivalent to $A$, and $k_\circ \in K(X)_{s+1}$ so that $k_\circ
A_\circ = A$.  Then $\Vert r_{\partial X} k_\circ \Vert_1 < c(R)$.
Suppose that $ r_{\partial X} k_\circ = \exp(\xi_\partial) $ for some
$\xi_\partial \in \k(\partial X)_1$, with $\xi_\partial$ taking values
in
\begin{equation} \label{alcm}
 \Alc^- := \{ \beta \in \Alc, \ \ \ \alpha_0(\beta) < 1 \}
 .\end{equation}
Then $\Vert \xi_\partial \Vert_0 < c$, hence $\Vert \xi_\partial
\Vert_1 < c(R).$
In fact, the image of $\Map(\partial X,\Alc^-)_1$ under the
exponential map is dense, since $K(\Alc - \Alc^-)$ is codimension $2$
in $K$.  It follows that $k_\partial = \exp(\xi_\partial)$, with
$\Vert \xi_\partial \Vert_1 < c(R)$.

Let $ \mE: \ \Omega^0(\partial X;\k)_1 \to \Omega^0(X;\k)_{\thh} $ be
any bounded linear operator with $r_{\partial X} \circ \mE = \on{Id}$ and
define $ \xi = \mE \xi_\partial , k_t = \exp(t \xi) .$ By Sobolev
multiplication $ \Vert k_t \Vert_{\thh} = \Vert \exp({t \xi})
\Vert_{\thh} \leq \exp(c \Vert t \xi \Vert_{\thh}) $ and similarly for
$k_t^{-1}$.  Define $ A_t = k_t A_0 .$ Then
$$ \Vert A_t \Vert_\hh = \Vert e^{\xi} A_0 \Vert_\hh \leq \Vert
e^{\xi} de^{-\xi} \Vert_\hh + \Vert \Ad(e^\xi) A_t \Vert_\hh < c(R) $$
and 
$$ \Vert \ddt A_t \Vert_\hh = \Vert \d_{A_t} \xi \Vert_\hh < c(R)
.$$ 
\end{proof}

The Proposition follows from Lemmas \ref{normlemma},\ref{pathlemma}
and \ref{surf}(b).

\section{The heat flow: existence of trajectories}

This section and the next one are modeled after R{\aa}de's treatment
\cite{ra:ym} of the Yang-Mills heat equation.  The norm-square of the
moment map
$$ f: \ \M(X)_s \to \R, \ \ [A] \mapsto \hh \Vert (r_+ - r_-)A
\Vert^2_{0} $$
is a $K(S)_{s + \hh}$-invariant smooth function.  We study the
equation
\begin{equation} \label{Aflow}
 \ddt [A_t] = -\grad(f)([A_t]) = - *_X d_A \mE_A \delta *_S (r_- -
 r_+) A .\end{equation}
Here $d_A$ is the covariant derivative for $A$, $\mE_A$ is the
harmonic extension operator for the generalized Laplacian $d_A^* d_A$,
$*_X, *_S$ are the Hodge star operators on $X, S$ respectively,
$\delta$ is the diagonal embedding, and the tangent space to $\M(X)_s$
is identified with $\ker d_A \oplus d_A^*$.

The gradient flow for $-f$ has the following description in local
slices.  For any $A,a \in \Omega^1(X;\k)_s$ there is a harmonic
extension operator
$$\mE_A^a: \ \Omega^0(\partial X;\k)_{s + \hh} \to \Omega^0({X},\k)_{s
+ 1} $$
such that
$$ r_{\partial X} \mE_A^a = \on{Id}, \ \ \ \ \d_A^*\d_{A+a}
\mE_A^a = 0 .$$
By definition of the local slices \eqref{UA} we have
$$T_{A + a} U_A = \ker \d_A^* \oplus \d_{A + a} .$$
Let $\pi_A^a$ be the orthogonal projection of $T_{A + a} U_{A+a}$
(with respect to the $L^2$ inner product) onto $T_{A + a} U_A$.
Define
$$ Q_A^a =  \pi_A^a \ \sx \, \d_{A+a} \, \mE_{A + a}^0 \delta *_S
(r_+ - r_-) .$$
The gradient flow equation is 
\begin{equation}
\label{Aflowinchart}
\ddt (A + a) = - Q_A^a(A+a), \ \ \ {A+a} \in U_A .
\end{equation}

\subsection{The linear initial value problem for the boundary data} 

Define linear operators
\begin{equation} \label{Ppm}
 P_{A,\pm}^a = \pm *_S r_\pm \pi_A^a \sx \d_{A+a} \mE^0_{A+a}
\delta \end{equation}
and set $P_{A,\pm} = P_{A,\pm}^0$.  The evolution equations for the
boundary data
\begin{equation} \label{Bpm}
 B_\pm := \frac{1}{2} *_S (r_{+} \pm r_{-}) A.
\end{equation} 
are
\begin{equation} \label{flowB}
\left(\frac{\d}{\d t} + P^a_{A,-} + P^a_{A,+} \right) B_- = 0, \ \ \
 \frac{\d}{\d t} B_+ = (P^a_{A,-} - P^a_{A,+}) B_-.
\end{equation}

\blem \label{Pops} For any smooth $A \in \A_\flat(X)$, the operators
${P}_{A,\pm}$ are elliptic pseudo-differential operators of order $1$
with principal symbol the same as the square root of the Laplacian on
$S$.  The sum ${P}_{A,+} + P_{A,-}$ is non-negative and self-adjoint
of order $1$.  \elem

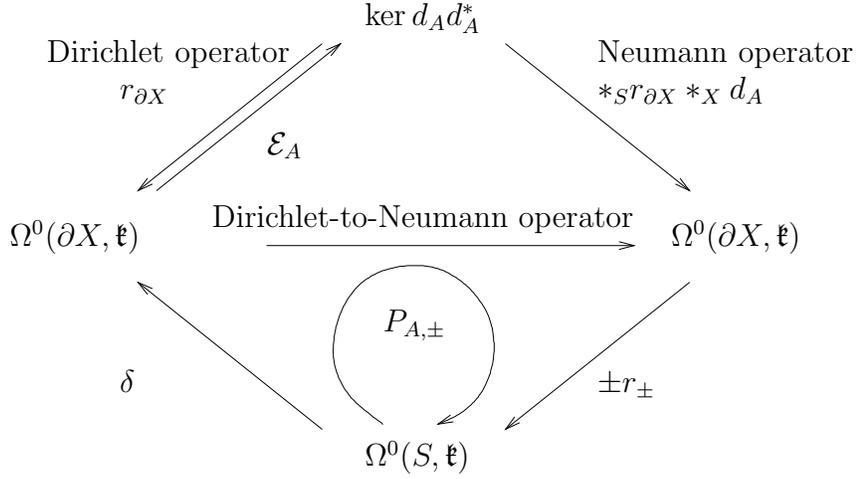
\begin{figure}[h]
\setlength{\unitlength}{0.00063333in}
\begingroup\makeatletter\ifx\SetFigFont\undefined%
\gdef\SetFigFont#1#2#3#4#5{%
  \reset@font\fontsize{#1}{#2pt}%
  \fontfamily{#3}\fontseries{#4}\fontshape{#5}%
  \selectfont}%
\fi\endgroup%
{\renewcommand{\dashlinestretch}{30}
\begin{picture}(7430,3829)(0,-10)
\path(2100,1858)(5100,1858)
\path(4980.000,1828.000)(5100.000,1858.000)(4980.000,1888.000)
\path(2550,3508)(1050,2308)
\path(1124.963,2406.389)(1050.000,2308.000)(1162.445,2359.537)
\path(4050,3508)(5550,2308)
\path(5437.555,2359.537)(5550.000,2308.000)(5475.037,2406.389)
\path(1162.445,1506.463)(1050.000,1558.000)(1124.963,1459.611)
\path(1050,1558)(2550,358)
\path(4124.963,456.389)(4050.000,358.000)(4162.445,409.537)
\path(4050,358)(5550,1558)
\path(2625.037,3409.611)(2700.000,3508.000)(2587.555,3456.463)
\path(2700,3508)(1200,2308)
\path(3047,398)(3044,400)(3038,405)
	(3027,413)(3011,425)(2992,440)
	(2970,457)(2947,475)(2924,494)
	(2901,512)(2881,529)(2862,545)
	(2845,560)(2829,574)(2815,588)
	(2802,601)(2790,614)(2779,627)
	(2768,641)(2758,655)(2747,669)
	(2737,684)(2727,700)(2718,717)
	(2708,734)(2700,752)(2691,770)
	(2683,788)(2676,807)(2670,825)
	(2664,844)(2659,862)(2655,880)
	(2651,898)(2648,915)(2646,933)
	(2644,951)(2643,969)(2642,988)
	(2642,1007)(2643,1027)(2644,1047)
	(2645,1068)(2648,1090)(2650,1111)
	(2654,1132)(2657,1153)(2662,1174)
	(2666,1194)(2671,1213)(2676,1232)
	(2682,1250)(2687,1267)(2693,1284)
	(2700,1302)(2708,1320)(2716,1339)
	(2725,1357)(2735,1375)(2745,1394)
	(2756,1412)(2767,1430)(2779,1447)
	(2791,1464)(2803,1480)(2815,1495)
	(2828,1509)(2840,1522)(2852,1534)
	(2865,1545)(2877,1556)(2890,1567)
	(2904,1577)(2919,1587)(2934,1597)
	(2950,1607)(2967,1616)(2984,1625)
	(3002,1634)(3020,1641)(3038,1649)
	(3056,1655)(3074,1661)(3091,1666)
	(3109,1671)(3127,1675)(3143,1679)
	(3159,1682)(3176,1685)(3194,1688)
	(3212,1690)(3232,1692)(3251,1693)
	(3271,1695)(3291,1695)(3312,1695)
	(3332,1695)(3351,1694)(3371,1693)
	(3389,1691)(3407,1689)(3425,1687)
	(3442,1684)(3458,1680)(3476,1675)
	(3495,1670)(3513,1664)(3531,1657)
	(3550,1649)(3569,1640)(3588,1631)
	(3607,1620)(3625,1609)(3643,1598)
	(3661,1586)(3677,1574)(3693,1562)
	(3708,1550)(3722,1538)(3736,1525)
	(3748,1514)(3760,1502)(3772,1489)
	(3785,1476)(3797,1462)(3809,1447)
	(3821,1432)(3833,1417)(3844,1400)
	(3855,1384)(3866,1367)(3876,1350)
	(3885,1333)(3894,1317)(3901,1300)
	(3909,1284)(3915,1267)(3921,1251)
	(3926,1234)(3931,1217)(3936,1199)
	(3940,1181)(3944,1162)(3947,1142)
	(3950,1122)(3952,1102)(3953,1081)
	(3955,1060)(3955,1039)(3955,1019)
	(3955,999)(3954,980)(3953,962)
	(3951,944)(3949,927)(3946,910)
	(3943,891)(3939,872)(3934,854)
	(3929,835)(3924,815)(3917,796)
	(3910,777)(3903,758)(3895,740)
	(3887,722)(3879,705)(3870,689)
	(3862,674)(3853,660)(3845,647)
	(3836,635)(3826,621)(3816,608)
	(3805,595)(3793,582)(3781,569)
	(3768,557)(3755,545)(3742,534)
	(3729,523)(3717,513)(3704,504)
	(3692,496)(3681,488)(3669,481)
	(3656,473)(3642,465)(3628,457)
	(3611,449)(3593,440)(3573,431)
	(3552,421)(3532,412)(3499,398)
\path(3597.753,472.483)(3499.000,398.000)(3621.186,417.249)
\put(1650,2008){\makebox(0,0)[lb]{\smash{{{Dirichlet-to-Neumann operator}}}}}
\put(2100,2608){\makebox(0,0)[lb]{\smash{{{$\mE_A$}}}}}
\put(300,3358){\makebox(0,0)[lb]{\smash{{{Dirichlet operator}}}}}
\put(4800,3358){\makebox(0,0)[lb]{\smash{{{Neumann operator}}}}}
\put(900,3058){\makebox(0,0)[lb]{\smash{{{$r_{\partial X}$}}}}}
\put(4800,3058){\makebox(0,0)[lb]{\smash{{{$*_S r_{\partial X} *_Xd_A$ }}}}}
\put(2900,3658){\makebox(0,0)[lb]{\smash{{{$\ker d_Ad_A^*$}}}}}
\put(5400,1858){\makebox(0,0)[lb]{\smash{{{$\Omega^0(\partial X,\k)$ }}}}}
\put(0,1858){\makebox(0,0)[lb]{\smash{{{$\Omega^0(\partial X,\k)$}}}}}
\put(2900,58){\makebox(0,0)[lb]{\smash{{{$\Omega^0(S,\k)$ }}}}}
\put(3050,1158){\makebox(0,0)[lb]{\smash{{{$P_{A,\pm}$}}}}}
\put(4800,658){\makebox(0,0)[lb]{\smash{{{$\pm r_\pm$}}}}}
\put(900,658){\makebox(0,0)[lb]{\smash{{{{$\delta$}}}}}}
\end{picture}
}
\caption{The operators $P_{A,\pm}$ \label{DNfig}}
\end{figure}

\begin{proof}  In the definition of
$P_{A,\pm}$, the Hodge star $\sx$ has the effect of exchanging tangent
and normal directions to the boundary.  The Dirichlet-to-Neumann
operator for the generalized Laplacian $d_A^* d_A$ is an elliptic
pseudo-differential operator of order one with principal symbol equal
to that of $d_A^* d_A$, see \cite[Chapter 21]{ho:an3},\cite{se:si}.
The operators $r_\pm$ and $\delta$ are Fourier integral operators of
order $0$, whose composition is the identity.  The composition
$P_{A,\pm}$ of the Dirichlet-to-Neumann operator with diagonal
embedding and restriction to $S_\pm$ is therefore a Fourier integral
operator with diagonal canonical relation, that is, a
pseudo-differential operator; see \cite{ho:an4}.  The relation between
the operators is shown in Figure \ref{DNfig}.  ${P}_{A,+} + P_{A,-}$
is non-negative:
\begin{eqnarray*}
\int_S \Tr( (P_{A_+} + P_{A,-})b \wedge *_S b) &=&
-\int_{\partial X}  \Tr( *_{\partial X} r_{\partial X} \sx \d_{A}
\mE_{A}^0 \delta b \wedge *_{\partial X}
\delta  b)  \\
&=& 
-\int_X d \Tr(\sx \d_{A} \mE_{A}^0 \delta b
\wedge \mE_A^0 \delta b) \\
&=& 
-\int_X \Tr(\sx \d_{A} \mE_{A}^0  \delta b \wedge 
\d_A \mE_{A}^0  \delta b) \\
& \geq&  0 .\end{eqnarray*}
The proof that $P_{A,+} + P_{A,-}$ is self-adjoint is similar.
\end{proof}

Since $P_{+,A}$ and $P_{-,A}$ have the same principal symbol, $P_{+,A}
+ P_{-,A}$ is elliptic (and therefore Fredholm) and $P_{+,A} -
P_{-,A}$ is a pseudo-differential operator of order $0$.  For $A$ not
smooth we will show the following properties of $P_{\pm,A}$.

\begin{proposition} \label{nonsmooth} For any $s > \hh$, 
flat $A \in \A(X)_s$ and $A + a \in U_A$,
\ben 
\item $P_{A,+}^a + P_{A,-}^a$ gives a Fredholm operator $\Omega^0(\partial
X;\k)_{s - \hh} \to \Omega^0(\partial X;\k)_{s-\thh}$;
\item $P_{A,-}^a - P_{A,-}^a$ gives a bounded operator
$\Omega^0(\partial X;\k)_{s - \hh} \to \Omega^0(\partial
X;\k)_{\min(2s-\thh,s-\hh)}$.  \een \end{proposition}
This will be derived from:
\blem \label{order} \label{vary} Let $s> \hh$.
\begin{enumerate}
\item Let $A_u$ be a differentiable path of flat connections.  Then
$\ddu \mE^0_{A_u} |_{u=v}$ is a bounded linear operator $
\Omega^0(\partial X;\k)_{s - \hh} \to \Omega^0(X;\k)_{\min(2s,s+1)}$;
\item Let $A_0 + a_u$ be a differentiable path of flat connections
such that $a_u$ lies in $U_{A_0}$. Then 
\ben
\item
$\ddu \d_{A_0 + a_u} |_{u =
v}$ is a bounded linear operator $ \Omega^0(X;\k)_{s + 1} \to
\Omega^1(X;\k)_{\min(2s,s+1)}$; 
\item $\ddu \pi_{A_0}^{a_u} |_{u = v}$ is a bounded linear operator $
(T_{A_u})_{s} \to \Omega^1(X;\k)_{\min(2s,s+1)}$,
\item $ \ddu P_{A_0,\pm}^{a_u} \alpha |_{u = v}$ is a bounded linear
operator $ \Omega^0(\partial X;\k)_{s - \hh} \to \Omega^0(\partial
X;\k)_{\min(2s-\thh,s -\hh)}$.  \een
\een \elem

\begin{proof} We have
\begin{eqnarray*}
 0 &=& \ddu \left( \d_{A_u}^* \d_{A_u} \mE_{A_u}^{0} b \right) \\ 
&=& \left( \sx \ad(\ddu {A_u}) \sx \d_{A_u} \mE_{A_u}^{0} +
\d_{A_u}^* \ad(\ddu {A_u}) \mE_{A_u}^{0} + \d_{A_u}^* \d_{A_u} \ddu \mE_{A_u}^0
\right) b.
\end{eqnarray*}
Hence
\begin{equation} \label{bigE} \d_{A_u}^* \d_{A_u} \ddu \mE_{A_u}^{0} b= - \left( \sx \ad(\ddu {A_u})
\sx \d_{A_u} + \d_{A_u}^* \ad(\ddu
(A_u)) \right) \mE_{A_u}^{0} b \end{equation} 
and
$$ 0 = r_{\partial} \left( \ddu \mE_{A_u}^{a_u} b \right) .$$
(a) now follows from elliptic regularity.  (b) (i) follows from $\d_{A_0 +
a_u} = \d_{A_0} + \ad(a_u)$ and Sobolev multiplication.  (ii) By
definition $\pi_{A_0}^{a_u} \alpha = \alpha - \d_{A_u}S \alpha $ where
$$S_u: \ (T_{A_u})_s \to (\ker \d_{A_u} \oplus \d_{A_0}^*)_s, \ \ \
\d_{A_0}^* (\alpha - \d_{A_u} S_u \alpha) = 0 . $$
Differentiating we obtain
$$ \d_{A_0}^* \ad(\ddu A_u) S_u \alpha = - \d_{A_0}^* \d_{A_u} \ddu
S_u \alpha .$$
We can re-write this 
$$ \d_{A_0}^* \d_{A_0} \ddu S_u \alpha = - \d_{A_0}^* (\ad(\ddu A_u) S_u
\alpha - \ad(A_u) \ddu S_u \alpha .$$
We apply elliptic regularity: Since the right-hand side is order
$2s-1$, $ \ddu S_u \alpha$ is order $2s + 1$, so $\ddu
\pi_{A_0}^{a_u}$ is order $2s$.
(iii) follows from (a), (b):(i),(ii).  
\end{proof}

Lemma \ref{nonsmooth} follows from \ref{vary} by choosing a path $A_u$
of connections from a smooth connection $A_0$ to $A$.  We can also use
\ref{vary} to derive bounds on the operator $\mE_A$.  

\blem For any $s > \hh$ and $\eps > 0$ there exists $\delta > 0 $ such
that if $A_u, \ u \in [0,1]$ is a differentiable path of flat
connections and
\begin{equation} \label{assume}
 \Vert A_u \Vert_{s} < \delta, \ \ \Vert \ddu A_u \Vert_{s} < \delta
\end{equation}
then $ \mE_{A_u}: \Omega^0(\partial X;\k)_{s + \hh} \to
\Omega^0(X,\k)_{s + 1}$ satisfies 
$$ \Vert \mE_{A_u} \Vert \leq \eps \Vert \mE_{A_0} \Vert.$$
Furthermore, $ P_{A_u,+} + P_{A_u,-}: \ \Omega^0(S;\k)_{s + \hh} \to
\Omega^0(S;\k)_{s - \hh}$ satisfies
$$\Vert P_{A_u,+} + P_{A_u,-} \Vert \leq \eps \Vert P_{A_0,+} +
P_{A_0,-} \Vert $$ 
and $ (P_{A_u,-} - P_{A_u,+}): \ \Omega^0(S;\k)_{s +
\hh} \to \Omega^0(S;\k)_{2s - \hh}$ satisfies 
$$\Vert P_{A_u,-} -
P_{A_u,+} \Vert \leq \eps \Vert P_{A_0,-} - P_{A_0,+} \Vert .$$ \elem

\begin{proof}  Apply \eqref{assume} to estimate the right-hand-side of 
\eqref{bigE}.  We obtain
$$ \Vert \ddu \mE_{A_u} b \Vert_s \leq \delta \Vert \mE_{A_u} b \Vert_s .$$
Hence
$$ \ddu \ln \Vert \mE_{A_u} b \Vert_s 
\leq \frac{ \Vert \ddu \mE_{A_u} b \Vert_s}{ \Vert \mE_{A_u} b
\Vert_s} \leq \delta .$$
Integrating with respect to $u$ gives
$$ \Vert \mE_{A_1} b \Vert_s \leq \Vert \mE_{A_0} b \Vert_s \exp(\delta).$$
The other estimates follows from the first, and the estimates on
$\d_{A_u}, r_{\partial X} $.
\end{proof}

From Lemma \ref{pathlemma} we get

\blem \label{PQbound} For any $s > \hh$, there exists a constant
$c(R)$ such that for any flat connection $A \in \Omega^1(X;\k)_s$ with
$f([A]) < R $, there exists a gauge transformation $k \in K(X)_{s +
1}$ such that $r_+ k = r_- k$ and
\begin{enumerate}
\item $ \Vert P_{+,k \cdot A} \Vert \leq c(R) $ as an operator of
order $-1$, and
\item  $\Vert P_{-,k \cdot A} \Vert \leq c(R) $ as
an operator of order $\min(0,s - 1)$. 
\end{enumerate}
 \elem

We solve the linear, time-independent boundary initial-value problem
\begin{equation}
\label{Bind+}
 ( \ddt + P_{A_0,+} + P_{A_0,-}) B_- = 0, \ \ \ \ 
\ddt B_+ = (P_{A_0,-} - P_{A_0,+}) B_-.
\end{equation}
We denote by $\Omega^1(S,\k)_{r,s}$ the Sobolev space of $\k$-valued
one-forms on $S$, time-dependent on an interval $[0,T]$, with $r$
derivatives in the time direction and $s$ derivatives on $X$; see the
appendix.  We assume throughout that $T < 1$.  By \ref{Sob} there is a
solution operator for \eqref{Bind+},
\begin{equation} \label{bsolve}
 M_{A_0,-}: \ \Omega^1(S,\k)_{s-\hh} \to \Omega^1(S,\k)_{(\hh - r, s -
   \hh + r)}, \ \ B_-(0) \mapsto B_- \end{equation}
for any real $r$, with norm $\Vert M_{A_0,-} \Vert < c(R) T^{r}$.
Suppose $\hh < s \leq 1$.  By Lemma \ref{order} $P_{A_0,-} -
P_{A_0,+}$ is an operator of order $s-1$ so that
$$ (P_{A_0,-} - P_{A_0,+}) B_- \in \Omega^1(S,\k)_{(\hh - r, 2s + r -
\thh )}.$$
Suppose now $r \in (0,1)$.  By \ref{Sob}(a) $\Omega^1(S,\k)_{(0,\hh +
  r, 2s - r - \thh )}$ is equal to $\Omega^1(S,\k)_{(\hh - r, 2s + r -
  \thh )}.$ By \ref{Sob}(g) integration gives a solution
$$B_+ \in \Omega^1(S,\k)_{0,\thh - r, 2s + r - \thh}.$$
to \eqref{Bind+}.  $\Omega^1(S,\k)_{(\thh - r, 2s + r - \thh)}$ embeds
into $\Omega^1(S,\k)_{(\thh - r, s + r - \thh)}$.  Taking $r = - \eps
$ and $r = 1 - \eps$ gives a solution
$$ B_- \in \Omega^1(S,\k)_{(\hh + \eps, s - \hh - \eps) \cap (-\hh +
\eps, s + \hh - \eps)}, \ \ \ B_+ \in \Omega^1(S,\k)_{\hh + \eps, 2s -
\hh - \eps}.$$

\subsection{Digression on gauge slices}

Later we will need slices for groups of gauge transformations of
various types.  Define
$$ T_A := \ker \d_{A}^* \oplus \d_{A}: \ \Omega^1(X;\k)_s \mapsto
(\Omega^0 \oplus \Omega^2)(X;\k)_{s-1} .$$
For $A$ smooth the trace map
$$ r_{\partial X}: T_A \to \Omega^1(\partial X;\k)_{s- \hh}
$$
is well-defined \cite[Chapter 13]{bo:el}.  

\blem \label{split} If $A \in \A_{\flat}(X)$ is smooth then there exist
$L^2$-orthogonal splittings
\begin{enumerate}
\item $T_A = (T_A \cap \ker (r_\partial  \sx) ) \oplus (T_A \cap \im
(d_A  \mE_A)  ) ;$
\item $ T_A = (T_A \cap \ker ((r_- - r_+)  \sx)) 
\oplus (T_A \cap  \im (d_A  \mE_A  \delta)) ;$
\item  $ T_A = (T_A \cap \ker (((r_- - r_+)  \sx) \times  (r_- - r_+)))
\oplus (T_A \cap \im ( ( d_A  \mE_A  \delta) \times ( \sx 
d_A  \mE_A  \delta)) .$
\end{enumerate}
\elem

\begin{proof}  We will prove (c); the others are similar.  
The reader may wish to compare this with the standard argument on
e.g. \cite[p. 194]{ho:an3}.  The adjoint of
\begin{equation} \label{op}
 ( d_A  \mE_A  \delta) \times ( \sx 
d_A  \mE_A  \delta): \ 
 \Omega^0(S;\k)_{s + \hh}^2 \to \ker(d_A \oplus d_A^*)_s, 
\end{equation}
is
\begin{equation} \label{adj}
((r_- - r_+)  \sx) \times (r_- - r_+)): \ \ker(d_A \oplus
d_A^*)_{- s} \to \Omega^1(S;\k)_{-s- \hh}^2 .\end{equation}
The pair \eqref{op},\eqref{adj} is an elliptic boundary problem,
(\eqref{adj} considered as a system of operators on $S$ is a
differential operator of order $0$) so if $a$ lies in the kernel of
\eqref{adj} with Sobolev class $s'$, then $a$ has Sobolev class
$\min(s' + 1)$ by Sobolev multiplication and elliptic regularity.
Starting with $s' = - s$ and repeating shows that $a$ is smooth.
Hence the kernel of \eqref{adj} for $s' =-s $ is identical to that for
$s' = s$, and $ \ker(d_A \oplus d_A^*)_s$ is the $L^2$-orthogonal sum
of the image of \eqref{op} and the kernel of \eqref{adj}.
\end{proof} 

\blem Let $A \in \A_{\flat}(X)$ be smooth.
\label{Coulomb}  
\ben
\item There exists a constant $\eps$ depending only on $\Vert
\d_A^{-1} \Vert, \Vert \d_A \Vert$ such that for any connection $A' \in
\A(X)_s$ satisfying $\Vert A' - A \Vert_s < \eps$, there exists a gauge
transformation $k \in K(X)_{s+1}$ such that $\d_A^* (k \cdot A' - A) =
0 $ and $r_{\partial X} \sx (k \cdot A' - A) = 0$.  
\item There exists a constant $\eps$ depending only on $\Vert
\d_A^{-1} \Vert, \Vert \d_A \Vert$ such that for any connection $A'
\in \A(X)_s$ satisfying $\Vert A' - A \Vert_s < \eps$, there exists a
gauge transformation $k \in K(X)_{s+1}$ such that $r_+ k = r_- k$ and
$k \cdot A' - A$ lies in
$$ \ker \d_A^* \oplus (r_- - r_+) \sx .$$
\item There exists a constant $\eps$ such that for any $A' \in
\A(X)_s$ satisfying $\Vert A' - A \Vert_s < \eps$, there exists a $g
\in G(X)_{s+1}$ with $r_+ g = r_-g$ depending smoothly on $A'$ such
that $g \cdot A' - A$ lies in
\begin{equation} \label{Cg}
\ker d_A \oplus d_A^* \oplus (r_- - r_+) \sx \oplus (r_- - r_+).
\end{equation}
\een \elem

\begin{proof}  By the implicit function theorem and Lemma
\ref{split}.
\end{proof}

\subsection{The linear initial-value problem in a local slice}

Let $A_0 \in \A_{\flat}(X)$ be smooth.  

\blem \label{perp} The image of $Q_{A_0}^0$ is perpendicular to
$T_{A_0} \cap \ker r_{\partial X}$.  \elem

\begin{proof}  Integration by parts gives
\begin{eqnarray*}
 \int_X \Tr(Q_{A_0}^0A \wedge \sx a) &=& \int_X \Tr(\sx \d_{A_0}
\mE_{A_0} *_S (\dr) A \wedge \sx a) \\ & = & \int_X \Tr(\d_{A_0}
\mE_{A_0} *_S (\dr) A \wedge a) \\ &=& \int_X \Tr(\mE_{A_0} *_S (r_- -
r_+) A \wedge \d_{A_0} a) \\ && + \int_{\partial X} \Tr(*_S (\dr) A
\wedge r_{\partial}a) \\ &=& 0 .
\end{eqnarray*}\end{proof}

We introduce norms on $T_{A_0}$
corresponding to the choice of different Sobolev norms in the
splitting
$$ (T_{A_0}) \cong \Ker r_{\partial X} |_{T_{A_0}} \oplus \Im (r_+ -
r_-) |_{T_{A_0}} \oplus \Im (r_- + r_+) |_{T_{A_0}} .$$
(To solve the Yang-Mills heat equation, one assumes $A \in H^s$ and
$F_A \in H^{s-1}$.)  The first summand is finite dimensional.  Define
\begin{multline} (T_{A_0})_{s}' = 
(\ker r_{\partial X})_{\hh + \eps,s}
\oplus 
(\Im (r_+ - r_-))_{(\hh + \eps, s - \hh - \eps) \cap (-\hh
+ \eps, s+ \hh - \eps)}
\\
\oplus 
(\Im (r_+ + r_-) )_{\hh + \eps, s - \hh - \eps} \end{multline}
where all operators are understood to be restricted to $T_{A_0}$.
Similarly, define
\begin{multline} (T_{A_0})_{0,s}' = ( \ker r_{\partial X})_{0,\hh + \eps,s}
\oplus 
(\Im ( r_+ - r_- ))_{0,(\hh + \eps, s - \hh - \eps)
\cap (-\hh + \eps, s+ \hh - \eps)} 
\\ \oplus 
(\Im (r_+ + r_- ))_{0,\hh + \eps, s - \hh - \eps} ;\end{multline}
\begin{multline} (T_{A_0})_{s}'' = ( \ker
r_{\partial X})_{-\hh + \eps} 
\oplus 
(\Im ( r_+ - r_-))_{ -\hh + \eps, s - \hh - \eps } 
\\ \oplus
(\Im ( r_+ + r_- ))_{- \hh + \eps, s - \hh - \eps } .\end{multline}
There is an embedding
$$ (T_{A_0})_{s}' \to (T_{A_0})_{\hh + \eps,s - \hh - \eps}, \ \ \
(b_0,b_-,b_+) \mapsto b_0 + \sigma(b_-,b_+);$$
there are similar embeddings for the other spaces.

\blem \label{Aindlemma} Solving the time-independent equation
\begin{equation} \label{Aind}
 (\ddt + Q_{A_0}^0) a = 0, \ \ a(0) = a_0 \end{equation}
defines an operator 
$$ M_{A_0}: \ (T_{A_0})_s \to \ (T_{A_0})'_s, \ \ \ a_0 \mapsto a $$
with $\Vert M_{A_0} \Vert < c(R) T^{-\eps}$.
\elem

\begin{proof} Let $A = \sigma(B_+,B_-)$ denote the unique lift of the solution
$(B_+,B_-)$ such that $A - A_0$ is perpendicular to $T_{A_0} \cap \ker
r_{\partial X}$.  By Lemma \ref{perp}, $A$ solves \eqref{Aind}.  The
estimate on the norm of $M_{A_0}$ follows from that on $M_{A_0,-}$ in
\eqref{bsolve}.
\end{proof}

\subsection{The non-linear initial value problem in a local slice}

\begin{theorem} \label{short}  For any Sobolev class $ s> \hh$
and smooth $A_0 \in \A(X)_{\flat}$, there exists a time $T$ depending
only on $R$ such that for sufficiently small $a_0 \in (T_{A_0})_{s} $,
the initial value problem \eqref{Aflowinchart} has a unique solution
$A = \varphi_{A_0}(a)$ on $[0,T]$ with $a$ in $ (T_{A_0})_s'$.  The
solution $A$ lies in $C^0([0,T],(U_{A_0})_s)$ and depends smoothly on
the initial condition $a_0$ in these topologies.
\end{theorem}

The proof uses standard Sobolev space techniques. For any flat
connection $A$ define a bounded linear operator
$$ L_{A}: (T_A)_{0,s}' \to (T_A)_{s}'' \ \ \ a \mapsto (\ddt + Q_{A}^0
)a. $$
That is,
$$ L_{A}(b_0,b_-,b_+) = ( \ddt b_0, (\ddt + P_{A,+} + P_{A,-})b_-,\ddt b_+
- (P_{A,-} - P_{A,+}) b_- ) .$$
Solving the inhomogeneous equation $(\ddt + P_{A_0,+} + P_{A_0,-}) u = f$
defines a right inverse to the operator
$$ \ddt + P_{A_0,+} + P_{A_0,-}: \ \Omega^1(S,\k)_{0,(\hh + \eps, s -
\eps) \cap (-\hh + \eps, s+ 1 - \eps)} \to \Omega^1(S,\k)_{-\hh +
\eps, s - \eps}, $$
with norm depending on $\Vert P_{A_0,+} + P_{A_0,-} \Vert$.  The
operators
$$ \ddt: \ \Omega^1(S,\k)_{0,\hh + \eps, s - \hh - \eps} \to
\Omega^1(S,\k)_{-\hh + \eps, s - \hh - \eps} $$
and
$$\ddt: \ (T_A \cap \ker r_{\partial X})_{\hh + \eps,s} \to (T_A \cap
\ker r_{\partial X})_{-\hh + \eps,s}
$$
are also invertible.  Therefore, $L_{A_0}$ has right inverse mapping
$(f_0,f_-,f_+)$ to
\begin{multline}
 ( (\ddt)^{-1}(f_0), (\ddt + P_{A,+} + P_{A,-})^{-1} (f_-), \\ 
(\ddt)^{-1}f_+ + (\ddt)^{-1} (P_{A,-} - P_{A,+}) (\ddt + P_{A,+} +
P_{A,-})^{-1} f_- ) .
\end{multline}
giving the solution to the time-independent inhomogeneous problem
\begin{equation} \label{inhom}
 (\ddt + Q_{A_0}^0 ) A = f, \ \ A(0) = 0. \end{equation}

\blem For $s \ge \hh$ and $A_0 \in \A(X)_s$ flat, there exists a gauge
transformation $k \in K(X)_{s+1}$ with $r_+ k = r_- k$ and $\Vert
L^{-1}_{k \cdot A_0} \Vert < c(R).$ \elem

\begin{proof}  This follows from Theorem \ref{PQbound}.
\end{proof}

The non-linear initial value problem is solved by perturbation.
Working in the slice $U_{A_0}$ near $A_0$ we write
$$A = \varphi_{A_0}(a_1 + a_2) $$
where $A_0 + a_1$ is the solution to the time-independent problem with
initial condition $a_1(0) = a_0$.  By Lemma \ref{chart}, the size of
$V_{A_0}$ is bounded from below by $c(R)$.  The problem
\eqref{Aflowinchart} becomes
$$ 0 = (\ddt + Q_A^a ) A = (\ddt + Q_{A_0}^0 - Q_{A_0}^0 +
 Q_{A_0}^{a_1 + a_2} ) (A_0 + a_1 + a_2) .$$
Since $A_0 + a_1$ solves the time-independent problem
\eqref{Aind}, 
\begin{equation} \label{A2}
( \ddt + Q_{A_0}^0) a_2 = (Q_{A_0}^0 - Q_{A_0}^{a_1 + a_2})(A_0 + a_1 + a_2) .\end{equation}
Define
$$ R_{A_0} = (Q_{A_0}^0 - Q_{A_0}^{a_1}) (A_0 + a_1) $$
and
$$ N_{A_0}a_2 = (Q_{A_0}^0 - Q_{A_0}^{a_1 + a_2}) (A_0 + a_1 + a_2) -
R_{A_0} .$$
$N_{A_0}$ is a non-linear operator with $N_{A_0} 0 = 0$.  We have to
solve the initial value problem
\begin{equation} \label{haveto}
(L_{A_0} - N_{A_0})a_2 = R_{A_0}, \ \ a_2(0) = 0 .\end{equation}
We will show that $N_{A_0}$ and $R_{A_0}$ have small norms for $T$
small.

\blem \label{Qvary} For $s > \hh$, the operators $\ddu Q_{A_0}^{a_u}
|_{u = v}$ and $Q_{A_0}^0 - Q_{A_0}^a $ have order $\min(0,s-1)$.
\elem

\begin{proof}  
(a) \bea
\ddu Q_{A_0}^{a_u} |_{u=v} &=& \ddu \left( \pi_{A_0}^{a_u} 
\sx \d_{A_u} \mE_{A_u} \delta \sx (r_+ - r_-)
\right) |_{u=v}  \\
&=& (\ddu  \pi_{A_0}^{a_u} |_{u=v}) \sx \d_{A_u} \mE_{A_u} 
\delta \sx (r_+ - r_-) \\
&+& \pi_{A_0}^{a_v} \sx \ad(a_v) \mE_{A_v} 
\delta \sx (r_+ - r_-)
\\
&+& 
\pi_{A_0}^{a_v} \sx \d_{A_v}  ( \ddu \mE_{A_u} |_{u=v}) 
(\pi_+ \times \pi_-)^* \sx (r_+ - r_-) .\eea
The result follows from Lemma \ref{vary}.  (b) Consider the path $A_0
+ ua, u \in [0,1]$, apply (a) and integrate with respect to $u$.
\end{proof}

\blem \label{depend} For $s \ge \hh$ and $T < c(R)$, the map $L_{A_0} -
N_{A_0}$ is a diffeomorphism of a neighborhood of $0$ in
$(T_{A_0})_{0,s}'$ onto a ball in $(T_{A_0})_{s}''$.  Furthermore,
there exists $k \in K(X)_{s + 1}$ with $r_+ k = r_- k $ such that
after replacing $A_0$ with $k\cdot A_0$ the radius of the ball is at
least $c(R) T^{-\hh + \eps}$ and the norm of $R_{A_0}$ is at most
$c(R)T^{\hh - 2\eps}$.  Equation \eqref{haveto} has a unique solution
$a_2$ with $\Vert a_2 \Vert''_s < c(R)T^{\eps - \hh}$.  The solution
$a_2$ depends smoothly on the initial condition $a_0$.  \elem

\begin{proof} Estimate for $R_{A_0}$:  
By interpolation $(r_- - r_+)(A_0 + a_1)$ lies in
$(T_{A_0})_{\eps,s-\eps}$.  By Lemma \ref{Qvary} $ R_{A_0}$ lies in
$(T_{A_0})_{\eps,2s-1-\eps}$.  Since $s \ge \hh$, this embeds in
$(T_{A_0})_{-\hh + \eps,s-\hh-\eps}$, and the norm of the embedding is
at most $cT^{\hh}$.  By Lemma \ref{Aindlemma} the norm of $a_1 =
M_{A_0} A_0$ is bounded by $c(R) T^{- \eps}$.  Hence 
$$ \Vert R_{A_0}
\Vert_{-\hh + \eps,s- \hh - \eps} \leq c(R) T^{\hh - \eps} + c(R)
T^{\hh - 2\eps} \leq c(R) T^{\hh - 2\eps} .$$

Now consider $N_{A_0}a_2$.  We have
$$ (D_{a_1 + a_2} N_{A_0})(a) = - \ddu Q_{A_0}^{a_1 + a_2 + u a}(A_0
+ a_1 + a_2) |_{u = 0} + (Q_{A_0}^0 - Q_{A_0}^{a_1 + a_2}) (a) .$$
The operators in this expression are again of order $\min(0,s-1)$, and
the same argument as for $R_{A_0}$ shows that
$$ \Vert D_{a_1 + a_2} N_{A_0} \Vert \leq c(R) (\Vert A_0 \Vert +
\Vert a_1 \Vert + \Vert a_2 \Vert) T^{\hh} <
c(R)(T^{\hh - \eps} + \Vert a_2 \Vert T^{\hh}) .$$
Therefore, for $T < c(R)$ and $\Vert a_2 \Vert \leq \qq c(R) T^{-
\hh}$ we have
$$ \Vert D_{a_1 + a_2} N_{A_0} \Vert \leq \hh \Vert L_{A_0}^{-1} \Vert .$$
It follows that $D_{a_1 + a_2}(L_{A_0} - N_{A_0})$ is invertible and
\begin{equation} \label{normD}
 \Vert D_{a_1 + a_2}(L_{A_0} - N_{A_0})^{-1} \Vert \leq 2 \Vert L_{A_0}^{-1} 
\Vert .\end{equation}
For any $a_0,a_1$ in $(T_{A_0})_s'$ let $a_t =
(1-t)a_0 + ta_1 .$ Then
\begin{eqnarray*}
  L_{A_0}^{-1}  ( (L_{A_0} - N_{A_0})a_1 -  (L_{A_0} - N_{A_0})a_0) 
& = &  (a_1 - a_0) \\ && - \int_0^1 L_{A_0}^{-1} (\d_{A_t} N_{A_0}) (a_1 -
  a_0) \d t \\
& \ge & \hh \Vert a_1 - a_0 \Vert .\end{eqnarray*}

We wish to show that the map $L_{A_0} - N_{A_0}$ is a diffeomorphism
of a neighborhood of zero onto a ball of radius $ c(R)T^{- \hh}$.  Let
$b$ lie in $(T_{A_0})_s''$ with $ \Vert b \Vert \leq c(R) T^{-\hh}$.
For all $s \in [0,1]$ consider the equation
$$ (L_{A_0} - N_{A_0} ) a_s = s b .$$
Solving this is equivalent to solving
$$ \ddth a_\theta = (D(L_{A_0} - N_{A_0})^{-1})b .$$
By \eqref{normD} this has a solution $a_s$
with 
$$\Vert a_s \Vert \leq 2 s \Vert L_{A_0}^{-1} \Vert
\Vert b \Vert .$$ 

It remains to show that $a_2$ depends smoothly on $a_0$.  This follows
from the implicit function theorem for Banach spaces applied to the
map $(a_0,a_2) \mapsto ( (L_{A_0} - N_{A_0})a_2, a_0)$.
\end{proof}

\blem The solution $a_1 + a_2$ we have constructed actually lies in $
(T_{A_0})_{\hh + \eps,s}$ and therefore by Sobolev embedding also in $
C^0([0,T],(T_{A_0})_s).$ \elem

\begin{proof}
$a_1 \in (T_{A_0})_{\hh + \eps,s - \eps}$ implies that
$$ R_{A_0} = (Q_{A_0}^0 - D \varphi_{A_0}^{-1} Q_{A_0}^{a_1}
\varphi_{A_0}) (A_0 + a_1) \in (T_{A_0})_{\hh + 2\eps,2s - 2\eps - 1}
.$$
Since $s \ge \hh$ this embeds into $(T_{A_0})_{\hh + 2\eps,s - 2\eps -
\hh}$.  Similar arguments show $N_{A_0}a_2$ lies in the same space.
Using that \eqref{flowB} is a parabolic system we obtain for $\eps <
\qq$
$$ a_1 + a_2 \in (T_{A_0})_{\hh + 2\eps,s} $$
and $a_1 + a_2$ depends smoothly on $a_0$ in this topology. 
\end{proof}

\subsection{Uniqueness and long-time existence}

\begin{theorem}  The initial value problem \eqref{Aflow} has a unique 
solution $[A] \in C^0_{\on{loc}}([0,\infty),\M(X)_s)$.
\end{theorem}

\begin{proof}  Since $\M(X)$ is dense in
$\M(X)_s$, we may after replacing $A_0$ with a gauge equivalent
connection choose a flat, smooth $A_0'$ arbitrarily close to $A_0$.
By Theorem \ref{short}, the solution to the heat flow in the slice at
$A_0'$ exists for a time $T$ depending only on an upper bound for
$f([A_t])$, which is non-increasing.  By iteration, a solution exists
for all times.  Let
$$ a, a' \in (T_{A_0})_s' \cap
C^0_{\on{loc}}([0,\infty), (T_{A_0})_s ) $$
be two solutions to \eqref{Aflowinchart}, with the same initial
connection $a_0 \in (T_{A_0})_s$.  Suppose that $a \neq a'$.  Let
$T_1$ be the largest number such that the restrictions of $a,a'$ to
$[0,T_1]$ are equal.  The restrictions of $a,a'$ to $[T_1,T]$ solve
the initial value problem \eqref{Aflowinchart} with initial data
$a(T_1) = a'(T_1) \in H^s$.  Without loss of generality we may assume
that $T_1 = 0$.  Let $a_1,a_1'$ denote the solutions to the
time-independent initial value problem \eqref{Aind}, and let $a_2 = a
- a_1$ and $a_2' = a' - a_1'$.  Since the solution to the
time-independent problem is unique, $a_1 = a_1'$.  For $\eps < \hh$,
the space $ (T_{A_0})_{0,s}' $ is the subspace of $(T_{A_0})_s'$ whose
elements vanish at $t = 0$, see \ref{Sob} (a).  It follows that $
a_2,a_2' \in (T_{A_0})_{0,s}' $ so that $a_2,a_2'$ solve \eqref{A2}.
The norms of the restrictions $A |_{[0,T_1]},A' |_{[0,T_1]} $ are
uniformly bounded as $T_1 \to 0$.  By the proof of existence, the
equation \eqref{A2} has a unique solution of norm less than $c(R)
T_1^{-\hh + \eps}$.  Therefore, for $T_1$ sufficiently small $a_2 =
a_2'$ on $[0,T_1]$, which is a contradiction.
\end{proof}

It would be interesting to know whether negative-time trajectories
exist, say for special $[A]$.  A natural candidate is those $[A]$
which extend to an open neighborhood of $X$, in a Riemann surface $X'$
containing $X$.

\section{The heat flow: Convergence at infinity}

The purpose of this section is to prove the following.

\begin{theorem}  \label{converge} For $s \ge 1$ and any $[A_0] \in
 \M(X)_s$, the trajectory $[A_t]$ converges in $\M(X)_s$ to a
critical point $[A_\infty]$ as $t \to \infty$.  For any critical
component $C$, the map 
$$ \rho_C: [0,\infty] \times \M(X)_C \to \M(X)_C, \ \ [A] \mapsto
[A_t] $$
is a deformation retract of $\M(X)_C$ onto $C$.
\end{theorem}
 
The critical points of $f$ are represented by flat connections $A$
such that
\begin{equation} \label{critf}
 \d_{A} \, \mE_{A}^0 *_S (r_{+} - r_{-}) A = 0 .\end{equation}

\blem \label{cancrit} For any $s \ge 1 $ and $[A] \in \crit(f)$, the
$K(S)$-orbit of any element $[A] \in \M(X)_{s,C}$ contains an element
$[A'] \in \crit(f)$ such that $A'$ is smooth and $r_\partial A'$ is a
harmonic, that is,
$$r_{-} A' = *_S \xi_-, \ \ \ \ r_+ A' = *_S \xi_+ $$
for some $\xi_\pm \in \t$.  The pair $(\xi_+,\xi_-)$ is uniquely
defined up to the action of $\Waff$ on $\t \oplus \t$.  \elem

\begin{proof}  
By Lemma \ref{circle} we may after replacing $A$ with a
gauge-equivalent connection assume $r_{+} A = *_S \xi_+$ for some
$\xi_+ \in \Alc .$ Since $[A]$ is infinitesimally fixed by $\xi$,
$r_{+} A$ and $ r_{-} A$ are also fixed, so $ \xi \in \k(S)_{r_+A} .$
$K(S)_{r_+A}$ is a compact Lie group, containing $T$ as a maximal
torus, and $\xi \in \k(S)_{r_+A}$.  It follows that there exists $k
\in K(S)_{r_+A}$ such that $k \cdot \xi \in \t$.  Then $k \cdot (r_-
A)$ is of the form $*_S \xi_-$, for some $\xi_- \in \t$.  The
intersection of the $K(S)$-orbit of $(r_- A,r_+ A)$ with $*_S\t \oplus
*_S \t$ is an orbit of $\Waff$.  It follows that $(\xi_+,\xi_-)$ is
unique up to the action of $\Waff$.  Smoothness of $A'$ follows
from elliptic regularity and bootstrapping.
\end{proof}

\blem \label{sequence} Let $[A_t]$ be a solution to the initial-value
problem \eqref{Aflow} in $\M(X)_s$, with $s \ge 1$.  There exists $t_i
\in \R$ and $k_i \in K(X)_{s+1}$ for $i = 1,2,\ldots$ such that $r_+
k_i = r_- k_i $ and
$$[k_i A_{t_i}] \to [A_\infty] $$
in $H^1$-norm, for some smooth $[A_\infty] \in \crit(f)$.  \elem

\begin{proof}  Since $ \ddt f(t) = \Vert Q_{A_t} A_t \Vert_0 $
and $f$ is bounded from below, there exists a sequence of times $t_i$
such that
$$ \ddt f(t_i) = \Vert Q_{A_{i}}A_{i} \Vert_0 \to 0 \ \text{as} \ 
i \to \infty.$$
where $A_i := A_{t_i}$.  By Lemma \ref{pathlemma}, there exists a
sequence of gauge transformations $k_i \in K(X)_{s+1}$ with
$r_+ k_i = r_- k_i $ such that
$$ r_+ (k_i \cdot A_i) \in *_S \Alc, \ \ \ \Vert k_i \cdot A_i
\Vert_{\hh} < c(R) \ \ \ \forall i =1,2,3,\ldots.$$
Choose $\eps \in (0,\hh)$.  Since the embedding $\A(X)_s \to \A(X)_{s
  -\eps}$ is compact, by passing to a subsequence and replacing $A_i$
with $k_i \cdot A_i$ we may assume that $A_i$ converges in $\A(X)_{s -
  \eps}$ to a flat connection $A_\infty \in \A(X)_{\hh - \eps}$.  Also
$r_+ (A_i)$ converges in any Sobolev norm to a connection
$A_{+,\infty}$.  By Stokes theorem
$$ \Vert Q_{A_{i}}A_{i} \Vert_0 = \int_{\partial X} \Tr( r_{\partial X}
*_X d_{A_i} \mE_{A_i}^0 B_{-,i} \wedge B_{-,i}).$$
Since $d_{A_i} \to d_{A_\infty}$ as operators of order $-1$,
$$ \int_{\partial X} \Tr( r_{\partial X} *_X d_{A_\infty} \mE_{A_\infty}^0
B_{-,i} \wedge B_{-,i}) \to 0 \ \text{as} \ i \to \infty .$$
Using Lemma \ref{Pops} we obtain
$$ \Vert r_{\partial X} *_X d_{A_\infty} \mE_{A_\infty}^0 *_S B_{-,i}
\Vert_{-\hh} \to 0 \ \text{as} \ i \to \infty .$$
Since $ r_{\partial X} *_X d_{A_\infty} \mE_{A_\infty}^0 *_S $ has
finite dimensional kernel and $\Vert B_{-,i} \Vert_0$ is bounded,
$B_{-,i}$ converges in $\Omega^1(S;\k)_{\hh}$ to some $B_{-,\infty}
\in \Omega^1(S;\k)_{\hh}$ and $r_-(A_i)$ converges to
$$ A_{-,\infty} = A_{+,\infty} + B_{-,\infty} \in
\Omega^1(S;\k)_{\hh} .$$  
Let $\ol{I}$ denote the closure of the set $I = \{ r_{\partial X} A_i
\}$.  Then $\ol{I}$ is compact (being the closure of the image of a
convergent sequence) and since $r_{\partial X}$ is proper $
r_{\partial X}^{-1} (\ol{I}) $ is compact.  Hence $[A_i]$ has a
subsequence converging to an element $[A_\infty]$ in $\M(X)_1$.
\end{proof}

We prove that the trajectory $A_t$ converges by showing that $
\int_0^\infty \Vert Q_{A_t} A_t \Vert_s \d t $ is finite.  Note that
$$ f(\infty) - f(0) = \int_0^\infty \ddt f = \int_0^\infty \Vert
Q_{A_t} A_t \Vert_0^2 $$
is finite.  We will give a lower bound for $\Vert Q_{A_t} A_t
\Vert_0$.  

\begin{theorem} \label{estlemma} Let $A \in \A(X)_\flat$ be a smooth connection
with $[A] \in \crit(f)$.  There exists $\gamma \in [\hh,1)$ such that
for any $A + a \in U_A$ sufficiently close to $A$ and $s \ge \hh$, 
\begin{equation} \label{est}
 \Vert Q_A^a(A+a) \Vert_{s - 1} \ge c | f(A+a) - f(A) |^\gamma .
\end{equation} 
\end{theorem}

\begin{proof}  The proof involves several lemmas.  Consider the
$L^2$-splitting
$$ T_{A+a} U_A = \ker (\d_{A+a} \oplus \d_A^* \oplus (r_- - r_+)\sx)
\oplus \Im \d_A \mE^a_A \delta ;$$
this is a variation of Lemma \ref{split} (c).  We denote the
projections by $\pi^0, \pi^1$ respectively.  Define
$$ \Sigma_A = U_A \cap \ker (r_- - r_+)\sx .$$
By Lemma \ref{Coulomb} (c) $\Sigma_A$ is a local slice for the $K(S)$-action
on $\M(X)_s$.  Consider the restriction $f|_{U_A^0}$, which has
$L^2$-gradient
$$ M: \ \Sigma_A \to T\Sigma_A, \ a \mapsto \pi^0 Q_A^a(A + a) .$$

\blem  If $\Vert a \Vert_s$ is small enough then
$ \Vert M(a) \Vert_{s-1} \ge \Vert Q_A^a(A + a) \Vert_{s-1} .$
\elem

\begin{proof} Note that $ \pi^1( Q_A^a(A + a)) = \d_{A+a} \mE_A^a \xi $
for some $\xi \in \k(S)$.  Also $\d_{A + a} = \d_A + \ad(a)$ implies
$\Vert \d_{A +a} \Vert \leq c \Vert \d_{A + a} \d_{A+a}^* \Vert$ for $
\Vert a \Vert_s < c'$.  So
\begin{eqnarray*}
 \Vert \pi^1 Q_A^a(A+a) \Vert_{s-1} &=& \Vert \d_{A+a} \mE_A^a \xi
 \Vert_{s-1} \\
&\leq& c \Vert \d_{A+a}^* \d_{A+a} \mE_A^a \xi \Vert_{s-2} \\
&\leq& c \Vert \d_{A+a}^* Q_A^a (A+a) \Vert_{s-2} \\
&\leq& c \Vert \d_A^* Q_A^a (A+a) \Vert_{s-2} + c \Vert \ad(a) \sx
Q_A^a (A+a) \Vert_{s-2} \\
&\leq& c \Vert a \Vert_s \Vert Q_A^a (A+a) \Vert_{s-1}.
\end{eqnarray*}
For $a$ sufficiently small we have
$$ c \Vert a \Vert_s \Vert Q_A^a (A+a) \Vert_{s-1} \leq \hh \Vert
Q_A^a (A+a) \Vert_{s-1}
$$
which proves the lemma. \end{proof}

It follows that $A + a$ is critical if and only if $M(a) = 0$.  The
derivative of $M(a)$ at $a = 0$ is the linear operator
$$ L: \ (\ker \d_A \oplus \d_A^* \oplus (r_- - r_+)\sx)_s \to
(\ker \d_A \oplus \d_A^* \oplus (r_- - r_+)\sx)_{s-1} $$
defined by
$$ L(\alpha) = \pi^0 \sx \left( \ad(\alpha) \mE^0_A (\dr) A + 
\d_A (\ddth \mE^{ \theta \alpha}_{A} ) (\dr) A + \d_A \mE^0_A
(\dr) \alpha \right).$$

\blem For any $\alpha \in \ker \d_A \oplus \d_A^* \oplus (r_+ -
r_-)\sx$,
$$ (\dr) \left( \ad(\alpha) \mE^0_A
 (\dr) A + \d_A (\ddth \mE^{\theta \alpha}_{A} ) (\dr) A +
 \d_A \mE^0_A (\dr) \alpha \right) = 0  .$$
\elem

\begin{proof}  We compute 
$$ (\dr) \ad(\alpha) \mE_A (\dr) A = 
\ad((\dr)\alpha) (\dr)A  $$
$$ (\dr) \d_A (\ddth \mE^{\theta \alpha}_{A} ) (\dr) A=
d_{r_+ A}  (\dr) A - d_{r_- A} (\dr) A = 0 $$
\begin{eqnarray*} (\dr) \d_A \mE^0_A (\dr) \alpha &=& 
d_{r_+ A} (\dr) \alpha  - d_{r_- A} (\dr) \alpha \\
&=& \ad( (\dr) A) (\dr) \alpha .
\end{eqnarray*}
\end{proof}

It follows that
$$ L(\alpha) = \sx \ad(\alpha) \mE^0_A (\dr) A + \sx \d_A (\ddth
\mE_{A}^{ \theta \alpha} ) (\dr) A + \sx \d_A \mE^0_A (\dr)
\alpha .$$

\blem $L$ is (1) self-adjoint and (2) Fredholm.
\elem

\begin{proof} (1) follows from 
\begin{eqnarray*}
 \int_X \Tr(\ad(\alpha_1) \mE^0_A (\dr) *_S A \wedge \alpha_2)
&=& \int_X \Tr(\mE^0_A (\dr) *_S A \wedge [\alpha_1,\alpha_2]) \\
&=& \int_X \Tr(\mE^0_A (\dr) *_S A \wedge [\alpha_2,\alpha_1]) \\
&=& \int_X \Tr(\ad(\alpha_2) \mE^0_A (\dr) *_S A \wedge \alpha_1) 
\end{eqnarray*}
by invariance of the inner product $\Tr$;
\begin{equation}
 \int_X \Tr(\d_A (\ddth \mE_{A}^{\theta \alpha_1} ) (\dr) A \wedge
\alpha_2) = \int_{\partial X} \Tr(r_{\partial X} \ddth \mE_{A}^{ 
\theta \alpha_1} (\dr) A \wedge r_{\partial X} \alpha_2) = 0
\end{equation}
using integration by parts and $r_{\partial X} \ddth
\mE_{A}^{\theta \alpha_1} = 0$; and 
\begin{eqnarray*}
 \int_X \Tr(\d_A \mE_A
*_S (\dr) \alpha_1  \wedge  \alpha_2)
&=& \int_{\partial X}\Tr( r_{\partial X} \mE_A *_S (\dr) \alpha_1 
\wedge r_{\partial X} \alpha_1)  \\
&=& \int_S \Tr(*_S (\dr) \alpha_1 \wedge (\dr) \alpha_2 )).
\end{eqnarray*}

(2) Let $\Sigma_A^0$ be the kernel of $L$, and $\Sigma_A^1$ its
$L^2$-orthogonal complement, so that $\Sigma_A = \Sigma_A^0 \oplus
\Sigma_A^1.$ The operator $r_- - r_+$ is Fredholm on $T_0\Sigma_A$, so
it suffices to show that $(r_- - r_+) L$ is Fredholm.  We have
\begin{multline}
(\dr)  L = (\dr) \sx ( \ad(\alpha) \mE^0_A (\dr)
A \\
+ \d_A (\ddth \mE^{ \theta \alpha}_{A} ) (\dr) A 
+ \d_A \mE^0_A
(\dr) \alpha ) \\
=  (\dr) \sx  ( \ad(\alpha) \mE^0_A (\dr)
*_S A \\
+ \d_A (\ddth \mE^{ \theta \alpha}_{A} ) (\dr) *_S A) +
P_{-,A}^0 (\dr) \alpha. 
\end{multline}
Since the class of Fredholm operators is closed under composition and
perturbation with compact operators, $ (r_- - r_+)  L$ is also
Fredholm.  This implies that $L$ is Fredholm.
\end{proof}

It follows that $\Sigma_A^0$ is finite dimensional and $L$ defines an
invertible operator $(\Sigma_A^1)_s \to (\Sigma_A^1)_{s-1}$.  Since
$L$ is the derivative of $M$ at $a=0$, it follows from the implicit
function theorem that there exists $\eps_1,\eps_2 > 0$ and a real
analytic map
$$ l: \ B_{\eps_1} \Sigma_A^0 \to B_{\eps_2} \Sigma_A^1 $$
such that $M(\alpha + l(\alpha)) = 0$.  Define
$$ f_0: \ B_\eps \Sigma_A^0 \to \R, \ \ \alpha \mapsto f(\alpha +
l(\alpha)) .$$
For any $\alpha_1 \in B_{\eps_1} \Sigma_A^0$ and $\alpha_2 \in \Sigma_A^0$ we
have
$$ ( \on{grad} f_0 (\alpha_1), \ \alpha_2)_0 = (M(\alpha_2 +
l(\alpha_2)),\alpha_2 + Dl(\alpha_1) \alpha_2)_0 .$$
Since $Dl(\alpha) \in \Sigma_A^1$, $ ( \grad f_0 (\alpha_1), \
\alpha_2)_{0} = (M(\alpha_2 + l(\alpha_2)),\alpha_2)_{0}$
which implies 
$$ \grad f_0(\alpha) = M(\alpha + l(\alpha)) .$$
Therefore the set of critical connections $M(a) = 0$ near $a = 0$ is
the set 
$$ \{ a = a_0 + l(a_0), a_0 \in B_{\eps_1} \Sigma_A^0, \ \grad f_0(a_0) = 0
\} .$$
For any $a$ sufficiently small, we may write $a = a_0 + l(a_0) + a_1$,
where
\begin{equation} \label{est2}
 \Vert a_0 \Vert_{\hh} \leq c \Vert a \Vert_{\hh}, \ \ \Vert l(a_0)
\Vert_{\hh} \leq c \Vert a \Vert_{\hh}, \ \ \Vert a_1 \Vert_{\hh} \leq c
\Vert a \Vert_{\hh} .\end{equation}
Now we estimate the left-hand side of \eqref{est}.  We have
\begin{eqnarray*}
 \pi^0 Q_A^a(A+a) &=& M(a) = M(a_0 + l(a_0) + a_1)\\
&=& \grad f_0 (a_0) + M(a_0 + l(a_0) + a_1)
- M(a_0 + l(a_0)) \\
&=& \grad f_0 (a_0) + \int_0^1 DM(a_0 + l(a_0) + 
s a_1 ) a_1 \d s \\
& =& \grad f_0(\delta_0) + L a_1 + L_1 a_1 
\end{eqnarray*}
where 
$$ L_1 = \int_0^1 (DM(a_0 + l(a_0) + s a_1) - DM(0)) a_1 \d s. $$
The spaces $\Sigma_A^0$ and $\Sigma_A^1$ are closed, disjoint subspaces of
$\Sigma_A$.  It follows that
$$ \Vert \pi^0 Q_A^a(A+a) \Vert_{s-1} \ge c \left(\Vert \grad
f_0(\delta_0) \Vert_{s-1} + \Vert L a_1 \Vert_{s-1} \right) - \Vert
L_1 a_1 \Vert_{s-1} .$$
From \eqref{est2} and the smooth dependence of $DM(a)$ on $a$ we see
that for $\Vert a \Vert_{s}$ sufficiently small
$$ \Vert \int_0^1 (DM(a_0 + l(a_0) + s a_1) - DM(0)) \d s \Vert \leq c
\eps_1 \Vert a_1 \Vert_{s} .$$
So for $\eps_1$ sufficiently small we have 
$$ \Vert \pi^0 Q_A^a (A+ a) \Vert_{s-1} \ge c \Vert \grad f_0 (a_0)
\Vert_{s-1} + c \Vert a_1 \Vert_{s} .$$
Since $M$ is the $L^2$-gradient of $E$,
\begin{eqnarray*} 
f(A+ a)  &=& f(A + a_0 + l(a_0) + a_1) \\
&=& f(A +a_0) + f(A + a_0 + l(a_0) + a_1) - f(a_0 + l(a_0)) \\
&=& f(A +a_0) + \int_0^1 (M(a_0 + l(a_0) + sa_1), a_1)_{0} \d s \\
&=& f(A +a_0) + (M(a_0 + l(a_0)),a_1)_{0} \\
&& + \int_0^1 \int_0^1 (DM(a_0 + l(a_0) + sta_1)sa_1,a_1)_{0} \d t \d s \\
&=& f(A + a_0) + (\grad f_0(a_0),a_1)_{0} + \hh(La_1,a_1)_{0}
+ (L_2 a_1,a_1)_{0},
\end{eqnarray*}
where
$$ L_2 = \int_0^1 \int_0^1 s(DM(a_0 + l(a_0) + sta_1) - DM(0)) \d t\d s .$$
The second term vanishes since $\nabla f_0(a_0) \in \Sigma_A^0, \ a_1 \in
\Sigma_A^1$.  The third term has norm at most $c \Vert a_1 \Vert^2_{\hh}$,
since $L$ is a bounded linear operator.  The norm of fourth term $
(L_2a_1,a_1)_{L_2}$ can be bounded in the same way as for $L_1$, by $c
\eps_1 \Vert a_1 \Vert^2_{\hh} .$ We conclude that for $\eps_1$
sufficiently small
$$ | f(A + a) - f(A) | \leq | f_0(a_0) | + c \Vert a_1 \Vert^2_{\hh} .$$
Since $f_0$ is real analytic, $\Sigma_A^0$ is finite-dimensional, $f_0(0) =
0$ and $\grad f_0 (0) = 0$ we conclude by a theorem of Lojasiewicz
\cite{lo:en} there exists $\gamma \in [\hh,1)$ such that for
sufficiently small $a_0$,
$$ \Vert \grad f_0(a_0) \Vert \ge c | f_0(a_0) |^\gamma .$$
(This is true for any norm on $\Sigma_A^0$.)  This completes the proof of
Theorem \ref{estlemma}.
\end{proof}

\blem \label{improve} Let $A_\infty$ be a representative of
$[A_\infty] \in \crit(f)$.  For any $\delta > 0$ and $T>\delta$, there
exists a constant $c$ such that if $A_t$ is a solution to the heat
equation \eqref{Aflowinchart} in the slice $U_{A_\infty}$, $0 < T_1
\leq T - \delta $ and $\Vert A_t - A_\infty \Vert_s \leq \eps$ for all
$ t\in [T_1,T]$ then
$$ \int_{T_1 + \delta}^T \Vert \ddt A \Vert_s \leq c \int_{T_1}^T
\Vert \ddt A \Vert_0 \d t .$$
\elem

\begin{proof}   Let $A'_t = \ddt A_t = - Q_{A_t}A_t$.  Then
$$ \ddt A'_t = - (\ddu Q_{A_u} )_{u=t} A_t + Q_{A_t} A'_t $$
so $A'_t$ satisfies the parabolic equation
$$ (\ddt - Q_{A_\infty}) A'_t =     - (\ddu Q_{A_u} )_{u=t} A_t 
+ (Q_{A_t} - Q_{A_\infty})A'_t .$$
In order to obtain estimates, we need to modify $A'_t$ to obtain a
function vanishing at $t = T_1$.  Let $\eta(t)$ be a smooth cut-off
function with $\eta = 0 $ on $[T_1,T_1+\delta/2]$ and $\eta = 1$ on
$[T_1 + \delta,T]$.  Then
$$ (\ddt - Q_{A_\infty}) (\eta A'_t) = (\ddu Q_{A_u} )_{u=t} (\eta A) +
(Q_{A_t} - Q_{A_\infty}) (\eta A'_t) + (\eta' A'_t) .$$
Hence $ \Vert \eta A'_t \Vert_{L^2([T_1,T],H^s)} $ is bounded by
\begin{multline*}
 c \Vert    - (\ddu Q^{A_u - A_\infty}_{A_\infty} )_{u=t} (\eta A) 
+ (Q_{A_t} - Q_{A_\infty})(\eta A'_t) + (\eta' A'_t) \Vert_{L^2([T_1,T],H^{
s-\hh})} \\
\leq c \Vert A'_t \Vert_{L^2([T_1,T],H^{s-\hh})} + c \Vert \eta A'_t
\Vert_{L^2([T_1,T],H^{ s-\hh})} + c \Vert \eta' A'_t
\Vert_{L^2([T_1,T],H^{ s-\hh})}
\end{multline*}
using \ref{order} and \ref{PQbound}. Using
H\"older's inequality we get 
\begin{eqnarray*}
 \Vert A'_t \Vert_{L^1([T_1+\delta,T],H^{s})}
&\leq& \Vert \eta A'_t \Vert_{L^1([T_1,T],H^{s})} \\
&\leq& (T-T_1)^\hh \Vert \eta A'_t \Vert_{L^2([T_1,T],H^{s})} \\
&\leq& c(R) (T-T_1)^\hh ( \Vert A'_t \Vert_{L^2([T_1,T],H^{s-\hh})} + \\
&& \Vert
\eta A'_t \Vert_{L^2([T_1,T],H^{s-\hh})} + c(R) \Vert \eta' A'_t
\Vert_{L^2([T_1,T],H^{s-\hh} )} ) \\
&\leq& c(R) (T-T_1)^\hh (1 + \delta^{-1}) \Vert A'_t
\Vert_{L^1([T_1,T],H^{s-\hh})} .\end{eqnarray*}
For $\hh > s > 0$ these norms are bounded from above by the
$L^2$-norm.  The case of arbitrary $s$ follows by bootstrapping.
\end{proof}

\blem \label{close} Let $[A_\infty] \in \crit(f)$.  Then there exist
constants $\eps_2 > 0, \hh \leq \gamma < 1, c > 0$ such that if
$[A_t]$ is a solution to the evolution equation \eqref{Aflow} and
$A_T$ is a representative of $[A_T]$ such that
$$ \Vert A_T - A_\infty \Vert_s \leq \eps_2 $$
for some $T> 0$, then either $f([A_t]) < f([A_\infty])$ for some $t >
T$ or $A_t$ is contained in the image of $U_{A_\infty}$ in $\M(X)$,
for all $t \ge T$ and $A_t$ converges in $(U_{A_\infty})_s$ to
$A_\infty'$ with $f(A_\infty') = f(A_\infty)$, as $t \to \infty$.  In
the second case,
\begin{equation} \label{closest}
\Vert A'_\infty - A_\infty \Vert_s \leq c \Vert A_T - A_\infty
\Vert_s.
\end{equation}
\elem

\begin{proof}  Assume that $f(A_t) > f(A_\infty)$ for all $t \in
[0,\infty)$.  Since $f$ is a smooth functional of $A$ and $A_\infty$
is a critical point, if we choose $\eps_2$ small enough then
\begin{equation} \label{fsmooth}
 | f(A_T) - f(A) | \leq c \Vert A_T - A_\infty \Vert^2_{s} .
\end{equation}
By Theorem \ref{short}, the solution to \eqref{Aflow} in
$C^0_{\on{loc}}([0,\infty),H^s)$ depends smoothly on the initial data
in $H^s$.  It follows that if $\eps_2$ is sufficiently small then
\begin{equation} \label{interval}
 \Vert A_t - A_\infty \Vert_{s} \leq c \Vert A_T - A_\infty \Vert_{s}
\end{equation}
for all $ t \in [T,T+1]$.  We claim that for $\eps_2 $ sufficiently
small, $\Vert A_t - A_\infty \Vert_{s} < \eps_1$ for all $t \ge T$.
Suppose the opposite. Let $T_1$ be the smallest number greater than
$T$ such that $ \Vert A_{T_1} - A_\infty \Vert_{s} \ge \eps_1 .$ By
\eqref{interval} if we choose $\eps_2$ small enough, then $T_1 > T +
1$.  Then by Lemma \ref{estlemma} for all $t \in [T,T_1]$ we have
\begin{eqnarray*}
\ddt (f(A_t) - f(A_\infty))^{1-\gamma}
&=& -(1-\gamma) \Vert Q_A^a(A+a) \Vert^2_0 \Vert (f(A_t) -
f(A_\infty))^{-\gamma} \\
&\leq& -c \Vert Q_A^a(A+a) \Vert_0 = -c \Vert \ddt A \Vert_0 .
\end{eqnarray*}
Integrating with respect to $t$ we get
\begin{eqnarray*}
 \int_T^{T_1} \Vert \ddt A \Vert_0 &\leq& c ( f(A_T) -
f(A_\infty))^{1- \gamma} \\
&\leq& c \Vert A_T - A_\infty \Vert_{s} \\
&\leq& c\eps_2^{2(1-\gamma)} \\
&\leq& c \eps_2 .\end{eqnarray*}
using \eqref{fsmooth}.  On the other hand,
\begin{eqnarray*}
\int_T^{T_1} \Vert \ddt A \Vert_{s} \d t 
& \geq & \Vert A_{T_1} - A_{T+1} \Vert_{s} \\
& \geq & \Vert A_{T_1} - A_\infty \Vert_{s} - \Vert A_{T+1} 
- A_\infty \Vert_{s} \\
&\ge& \eps_1 - c \eps_2 .\end{eqnarray*}
It follows from \ref{improve} that $\eps_1 - c\eps_2 \leq c
\eps_2$.  For $\eps_2$ sufficiently small, this gives a contradiction.

We conclude that for $\eps_2$ sufficiently small, $\Vert A_t -
A_\infty \Vert_{s} < \eps_1$ for all $t \ge T$.  Then
$$ \int_{T+1}^\infty \Vert \ddt A \Vert_{s} \leq c \int_T^\infty
\Vert \ddt A \Vert_{0} \d t \leq c (f(A_T) - f(A_\infty))^{1-\gamma} .$$
It follows that $A_t$ converges to $A_\infty'$ as $t \to \infty$.  By
Lemma \ref{cancrit}, the set of critical values of $f$ is locally
finite, so $f([A_\infty']) = f([A_\infty])$ for $\eps_2$ sufficiently
small.

It remains to prove the estimate \eqref{closest}.  This follows
from the computation
\begin{eqnarray*}
\Vert A_\infty' -  A_{T+1} \Vert_{s}
&\leq& \int_{T+1}^\infty \Vert \ddt A \Vert_{\hh} \\
&\leq& c(f(A_T) - f(A_\infty))^{1-\gamma} \\
&\leq& c \Vert  A_T-  A_\infty \Vert_{\hh}^{2(1-\gamma)} \\
&\leq& c \Vert A_T - A_\infty \Vert_{\hh}.
\end{eqnarray*}
\end{proof}

Proof of Theorem \ref{converge}. Let $t_i,k_i,[A_\infty]$ be a
sequence given by Proposition \ref{sequence}.  Since the equation
\eqref{Aflow} is invariant under $K(S)$, the trajectory $k_n [A_{t_n +
t}]$ is also a solution.  For $n$ sufficiently large, $k_n [A_{t_n +
t}]$ satisfies the assumptions of Proposition \ref{close}.  Therefore,
$k_n [A_{t_n + t}]$ converges to some $[A_\infty']$.  It follows that
$A_t \to k_n^{-1}[A_\infty']$.

It remains to show that $[A_\infty]$ depends continuously on the
initial data.  Let $\eps_1 > 0$.  Let $[A_t]$ be a solution to
\eqref{Aflow}, and $A_\infty$ a representative for $[A_\infty]$.  By
Proposition \ref{close}, there exists an $\eps_2 > 0$ such that if
$[A'_t]$ is another solution to \eqref{Aflow} such that $\Vert A'(T)-
A_\infty \Vert_s \leq \eps_2$ for some $T \ge 0$ and representative
$A'(T)$, and $f([A'_\infty]) = f([A_\infty])$, then $\Vert A'_\infty -
A_\infty \Vert_s \leq \eps_1$.  Choose $T$ sufficiently large so that
$\Vert A_T - A_\infty \Vert_s \leq \eps_2/2$, where $A_T$ is the
representative in the slice $U_{A_\infty}$.  By the first part of the
theorem, there exists $\eps_3 > 0$ such that if $\Vert A'_0 - A_0
\Vert_s < \eps_3$ then $\Vert A'_T - A_T \Vert_s < \eps_2/2$.  We
conclude for any $\eps_3>0 $ there exists $\eps_2 >0$ such that if
$\Vert A'_0 - A_0 \Vert_s < \eps_2$ and $f([A'_\infty]) =
f([A_\infty])$ then $\Vert A'_\infty - A_\infty \Vert_s < \eps_3$, for
some representatives $A'_\infty,A_\infty$.  This completes the proof.

\section{The stratification defined by the heat flow}

\subsection{The critical set for $f$}

Let $\cC$ be the set of connected components of $\crit(f)$.  For any
$C \in\cC$, define
$$\M(X)_C = \{ [A] \in \M(X), \ \ [A_\infty] \in C \} $$
so that
$$ \M(X) = \bigcup_{C \in\cC} \M(X)_C .$$

\blem For any $R > 0$, there are a finite number of critical sets $C$
whose image under $\Phi$ intersects a ball $B_R$ of radius $R$ in
$\Om^1(X;\k)_{s - \hh}$.  \elem

\begin{proof}  For any subgroup $H \subset T$,
let $F$ be the fixed point set of the action of $H$ on $\Phi^{-1}(B_R)
\cap r_{\partial X}^{-1}(*_S \t \oplus *_S \Alc)$.  By the symplectic
cross-section theorem for loop group actions \cite{me:lo} there are at
most a finite number of components of $F$.  Each contains at most one
component of the critical set of $\zeta$, given by the set of points
$[A] \in F$ such that $*_S \Phi(A)$ is perpendicular to the Lie
algebra $\h$.  By Lemma \ref{cancrit}, any critical component of $\zeta$
contains elements of this form.  Hence, the number is finite.
\end{proof}  

We will give an explicit description of $\crit(f)$.  By \eqref{critf},
any $[A] \in \crit(f)$ is fixed by a one-parameter subgroup, which we
may assume is generated by an element
$$\xi \in \t, \ \ \ *_S \xi = (r_- - r_+)A  .$$ 
The centralizer $K_\xi \subset K$ is a connected subgroup of $K$, with
Lie algebra $\k_\xi$.  Let $\M(X;K_\xi)$ denote the moduli space of
flat $\k_\xi$-connections on $X \times K$, modulo restricted gauge
transformations in $K_\xi$.  The inclusion $\k_\xi \to \k$ induces an
embedding 
$$ \iota_\xi: \ \M(X;K_\xi) \to \M(X) $$ 
whose image is the fixed point set $\M(X,K)^{U(1)_\xi}$ of $U(1)_\xi$.
Let
$$ \M(X;K_\xi;\xi) = \{ [A] \in \M(X;K_\xi), \ \ (r_- - r_+)A = *_S
\xi \}.$$
Then
\begin{equation} \label{critun}
 \crit(f) = \bigcup_{\xi} K(S) \cdot \iota_\xi \M(X;K_\xi;\xi) .
\end{equation} 

The quotient of $K_\xi(S) \backslash \M(X;K_\xi;\xi)$ by the loop
group $K_\xi(S)$ is homeomorphic to the moduli space
$\M(\ol{X};K_\xi;\xi)$ of bundles with constant central curvature
$\xi$.  The homeomorphism is constructed by subtracting a
$Z(K_\xi)$-connection $A_\xi$ with $(r_- - r_+)A_\xi = *_S \xi$; such
a central connection exists after adding a boundary component with
marking $\ol{\xi}$, the reflection of $\xi$ into the fundamental
alcove.  This gives a homeomorphism to the moduli space of flat
bundles on $\ol{X}$, with the additional marking.  A similar
discussion gives a homeomorphism between $ \M(\ol{X};K_\xi;\xi)$ and
the same space.  From \eqref{critun} we obtain
\begin{equation} \label{comparecrit}
 K(S) \backslash \crit(f) \to \bigcup_{\xi} \M(\ol{X};K_\xi;\xi)
\end{equation}
where the sum is over $\xi \in \t_+$ such that $\xi$ is a coweight for
$K_\xi/[K_\xi,K_\xi]$.  The same description holds for the critical
set of the Yang-Mills functional on $\A(\ol{X})$ \cite[Section
6]{at:mo}.

\subsection{The semistable stratum}

From now on we drop the Sobolev subscripts $_s$; the statements that
follow hold for any $s \ge 1$.  For the following compare
\cite[p.36]{mu:ge}.

\begin{definition}
The {\em semistable} locus $\M(X)^{\ss}$ is the set of all points
$[A]$ with $[A_\infty] \in \Phinv(0)$.  The {\em stable} locus
resp. $\M(X)^{\s}$ is the set of all points $[A] \in \M(X)$ with
$[A_\infty] \in \Phinv(0)$, $[A_\infty] \in G(S)[A]$, and $[A]$ has
finite stabilizer.
\end{definition}

The aim of this section is to prove the following theorem.

\begin{theorem} \label{semistable}  
\ben
\item $\M(X)^{\s}$ is open.
\item The map $\M(X)^{\s} \to K(S) \backslash (\Phinv(0) \cap
\M(X)^{\s})$ is $G(S)$-equivariant.
\item $\M(X)^{\s}$ is $G(S)$-invariant; hence $\M(X)^{\s}$ is the
union of $G(S)$-orbits in $\M(X)$ that intersect $\Phinv(0)$ and have
finite stabilizers.
\item $\M(X)^{\s}$ is the set of $[A] \in \M(X)^{\ss}$ such that $[A]$
has finite stabilizer and $G(S)[A]$ is closed in $\M(X)^{\ss}$.
\item If $\M(X)^{\s}$ is non-empty then $\M(X)^{\s}$
is dense in $\M(X)^{\s}$;
\item $\M(X)^{\ss}$ is a $G(S)$-invariant
open subset of $\M(X)$; and
%
%
\item the closures of orbits $G(S)[A_1],G(S)[A_2]$ intersect in
$\M(X)^{\ss}$ if and only if $K(S)[A_{1,\infty}] =
K(S)[A_{2,\infty}]$.
\item If $g \ge 2$, or if $g = 1$ and there is at least one generic
marking, or if $g = 0$ and there are at least three generic markings,
then $\M(X)^s$ is non-empty.  \een
\end{theorem}

\begin{remark}  Probably the assumption in (e)-(g) is unnecessary.
\end{remark}

The theorem depends on the convexity of a certain functional, which
was introduced in Guillemin-Sternberg \cite{gu:ge} and used by
Donaldson \cite{do:ne} in infinite dimensions.  For any $[A] \in
\M(X)$, choose an element $l $ in the fiber $\L(X)_{[A]}$ of the
pre-quantum line bundle over $\M(X)$.  The central $\C^*$ in
$\widehat{G(S)}$ acts trivially on $\L(X)$, and so the action induces
an action of $G(S)$.  Consider the smooth function
$${\Psi}: \ {G(S)} \to \R, \ \ \ \ \  g \mapsto - \ln \Vert g l \Vert^2 .$$
Since ${K(S)}$ acts on $\L(X)$ preserving the metric, ${\Psi}$
descends to a smooth function on ${G(S)}/{K(S)}$, which is diffeomorphic
to $\k(S)$ via the map
$$ \k(S) \to G(S)/K(S), \ \ \ \ \ \xi \mapsto [\exp(i\xi)]. $$ 
We denote by $\psi$ the induced function $\k(S) \to \R$.  As in the
finite-dimensional case the function $\psi$ satisfies
\begin{equation} \label{diff}
 \ppxi \psi = \Phi^\xi, \ \ (\ppxi)^2 \psi = \Vert \xi_{\M(X)} \Vert^2
\end{equation}
where $\Phi$ is to be evaluated at $\exp(i\nu)[A]$.  $\psi$ has
the following properties:

\blem \label{min} \ben
\item $\psi$ is convex.
\item The only critical points of $\psi$ are zeros of the
pull-back of $\Phi$.  
\item The gradient flow lines from $\psi$ map onto the gradient
flows lines for $f$ under the map $\xi \mapsto \exp(i\xi)[A]$.
\item If $[A_\infty] \notin \Phinv(0)$ then $\psi([A_t]) \to -\infty $
as $ t \to \infty$.
\een
\elem

The proofs of (a)-(c) are the same as in the finite-dimensional case.
It follows from (a) that if $[A] \in \Phinv(0)$ and $g \in G(S)$ is
such that $g [A] \in \Phinv(0)$, then $g \in K(S)$, so $ \k(S) \to
\M(X), \ \xi \mapsto \exp(i \xi) [A] $ is an immersion of Banach
manifolds.  Proof of (d): If $\xi_\infty = *_S\Phi([A_\infty])$, then
$$ \ddt \psi([A_t]) \to - \Phi([A_\infty])^{\xi_\infty} = - \Vert
\xi_\infty \Vert^2 \ \text{as} \ \ t \to \infty $$
which implies $\psi([A_t]) \to -\infty$.

\vskip .1in

Proof of Theorem \ref{semistable}: (a) $D \Phi$ is surjective at
$\Phinv(0) \cap \M(X)^{\s}$; it follows that $G(S) (\Phinv(0) \cap
\M(X)^{\s} )$ contains an open neighborhood of $\Phinv(0) \cap
\M(X)^{\s}$.  Therefore, $[A] \in \M(X)^{\s}$ if $[A_\infty] \in
\Phinv(0) \cap \M(X)^{\s}$.  By continuity of the retraction
$\M(X)^{\ss} \to \Phinv(0)$, $[A'_\infty]$ lies in $ \Phinv(0) \cap
\M(X)^{\s}$ for $[A']$ sufficiently close to $[A]$.  This implies
$[A'] \in \M(X)^{\s}$.

(b) Let $[A'] = g[A]$ for some $g \in G(X)$ and $[A] \in
\M(X)^{\s}$.  Convexity of $\psi$ implies that $\psi$ has at most one
critical point $\xi_\infty$ on $\k(S)$, and $\xi_\infty$ is a global
minimum for $\psi$.  By Lemma \ref{min} (b), $\exp(\xi_\infty[A]) \in
\Phinv(0)$.  By Lemma \ref{min} (d), $[A_\infty'] \in \Phinv(0)$.
Since the retraction $\M(X)^{\ss} \to \Phinv(0)$ is continuous,
$G(S)[A]$ is connected and lies in $\M(X)^{\ss}$, and $K(S)[A_\infty]$
is isolated in $\Phinv(0)$, we must have $[A'_\infty] \in
K(S)[A_\infty]$.

(c) follows immediately from (b).

(d) Suppose $g_t[A] \to [A']$ as $t \to \infty$, for some $[A] \in
\M(X)^{\s}$ and $[A'] \in \M(X)^{\ss}$.  By continuity of the
retraction, $[A_\infty] = k[A'_\infty]$, for some $k \in K(S)$.  This
implies $[A'] \in G(S)[A]$.  Conversely, if $[A] \in \M(X)^{\ss}$ has
finite stabilizer and $G(S)[A]$ is closed then $[A_\infty] \in
G(S)[A]$, which implies $[A] \in \M(X)^{\s}$.

(e) Suppose $[A] \in \M(X)^{\ss}$ and $S_\infty$ be the gauge slice
\eqref{Cg}.  It suffices to show that the intersection $\M(X)^{\s}
\cap S_\infty$ is dense in $S_\infty$.  In order for $G(S)[A']$ to be
closed, it suffices that the intersection $G(S)[A'] \cap S_\infty$ is
closed, again by continuity of the retraction.  Note $S_\infty$ is
$K(S)_{[A_\infty]}$-equivariantly isomorphic to the representation
$TS_\infty$, and the intersection $G(S)[A'] \cap S_\infty$ maps onto a
$G(S)_{[A_\infty]}$-orbit in $S_\infty$.  The union of closed
$G(S)_{[A_\infty]}$ orbits in $TS_\infty$ with finite stabilizer is
either empty or dense; it follows that $\M(X)^{\s}$ is either empty or
dense.

(f) Suppose that $[A] \in \M(X)^{\ss}$ and $g \in G(S)$.  Let $[A_i]$
be a sequence in $\M(X)^{\s}$ with $\lim([A_i]) = [A]$.  By part (a),
$g[A_i] \in \M(X)^{\s}$, which implies $g[A] \in \M(X)^{\ss}$ by
smooth dependence of $[A_\infty]$ on $[A]$.

(g) follows since the retraction $\rho_0: \ \M(X)^{\ss} \to \Phinv(0)$
is continuous, and the orbits of $K(S)$ in $\Phinv(0)$ are closed, see
\ref{surf}.

(h) is a well-known consequence e.g. of the holonomy description and
left to the reader.

\subsection{The unstable strata}

Let $C \in \cC$.  The purpose of this section is to prove

\begin{theorem} \label{main}
If $s \ge 1$ and the genus of $\ol{X}$ is at least $2$, then
$\M(X)_{C}$ is a smooth K\"ahler $G(S)_{s + \hh}$-invariant
submanifold of finite codimension.
\end{theorem}

Let $\xi_\pm = *_S r_\pm([A])$ as described in Lemma \ref{cancrit} and
$$\xi_C = \xi_+- \xi_- .$$
 Let
$$\C^*_C = \{ \exp(  \tau \xi_C), \ \ \tau \in \C \} $$
denote the one-parameter subgroup of $G(S)$ generated by $\xi_C$.  Let
$$Z_C \subseteq \M(X)^{\C^*_C} $$ 
be the component of the fixed point set of $\C^*_S$ containing $[A']$.
Let $P_C$ be the parabolic subgroup corresponding to the element
$\xi_C$, so that
$$ P_C = \{ g \in G, \ \ \lim_{\tau \to -\infty} \exp(i\tau \xi_C) g
\exp(- i\tau \xi_C ) \ \text{exists} \} .$$
Let $ P_C = L_C U_C $ be its standard Levi decomposition, so that
$L_C$ is the centralizer of $\xi_C$.  Let $K_C$ denote the maximal
compact subgroup of $L_C$, that is, $K_C = L_C \cap K .$ Let
$P_C(S),L_C(S),K_C(S)$ etc. denote the identity components of the loop
groups of maps $S \to P_C, L_C, K_C$ etc.  Let $\pi_C: \ \k \to \k_C$
denote the projection, and $\pi_C(S)$ the pointwise projection
$\pi_C(S): \ \k(S) \to \k_C(S)$.  The group $L_C(S)$ acts on $Z_C$,
and the action of the subgroup $K_C(S)$ is symplectic, with moment map
$\pi_C \Phi$.  The same argument as in the proof of Theorem
\ref{semistable} part (a) shows that if the set $Z_C^{\s}$ is
non-empty, then $Z_C^{\ss}$ is $L_C(S)$-invariant; this holds with the
same assumptions on the genus and markings as for the action of
$G(S)$.  Define
$$ Y_C = \{ [A] \in \M(X), \ \ \ g[A] \to Z_C \ \text{as} \ g \to 0 \
\text{in} \ \C^*_C \} .$$
\blem 
\ben
\item $Y_C$ is a $P_C(S)$-invariant complex submanifold.
\item  If $[A] \in G(S)Y_C$ then $f([A]) \ge f(C)$.
\een
\elem

\begin{proof} (a) $Y_C$ is a stable manifold for the gradient flow of 
$-(\Phi,\xi_C)$.  Since $(\Phi,\xi_C)$ is Morse-Bott, $Y_C$ is an
embedded complex submanifold by the stable manifold theorem
\cite[Theorem III.8]{sh:gl}.  To show $Y_C$ is invariant under
$P_C(C)$, suppose $ \exp(z\xi_C) y \to z$ for some $z \in Z_C$.  Then
$$ \exp(z \xi_C) py = \exp(z \xi_C) p \exp(-z \xi_C) \exp(z
\xi_C) y \to lz \in Z_C $$ 
since $L_C(S)$ commutes with $\xi_C$.  (b) Since $G(S) = K(S) P_C(S)$
and $Y_C$ is $P_C(S)$-invariant, $G(S)Y_C = K(S)Y_C$.  Since $f$ is
$K(S)$-invariant we may assume $ [A] \in Y_C$.  Then $ (\Phi(A),\xi_C)
\ge (\xi_C,\xi_C) $ which implies (b).
\end{proof}

Let $Y_C^{\ss}$ denote the $P_C(S)$-invariant set of $[A]$ in $Y_C$
such that $[A_\infty]$ lies in $Z_C^{\ss}$.

\blem \label{set} 
 If $[A] \in G(S)Y_C^{\ss}$ then (i) $\xi_C$ is the unique
closest point to $0$ in $\Phi(\ol{P_C(S)[A]})$ and (ii) $K(S)\xi_C$ is
the set of points closest to $0$ in $\Phi(\ol{G(S)[A]})$. 
\elem

\begin{proof} (i)
We have to show that $\xi_C$ lies in $\Phi(\ol{P_C(S)[A]})$.  Suppose
$\exp(it \xi_C) [A] \to [A']$ as $t \to -\infty$.  Then
$\ol{P_C(S)[A]}$ contains $\ol{L_C(S)[A']}$, so it suffices to show
that $\Phi(\ol{L_C(S)[A']})$ contains $\xi_C$.  But this is the
definition of $Z_C^{\ss}$.  The argument for (ii) is similar.
\end{proof}

\blem \label{open} $Y_C^{\ss}$ is the unique open $P_C(S)$-invariant
neighborhood of $Z_C \cap \Phinv(C)$ in $Y_C$.
\elem

\begin{proof}  
The closure of any $P_C(S)$-orbit in $Y_C^{\ss}$ intersects $\Phinv(C)
\cap Z_C$.
\end{proof}

\blem If $\nu \in \k(S)$ and $\nu_{[A]}$ is tangent to $T_{[A]} Y_C$
for some $[A] \in C$, then $\nu \in \k_C(S)$.  \elem

\begin{proof} We follow Kirwan
\cite[p.50]{ki:co}. We have
$$ \Phi( \exp(t\nu) [A]) = \xi_C + t[\nu,\xi_C] + e(t) $$
where $e(t) = O(t^2)$ as $t \to 0$.  This implies that 
$$ ( \Phi( \exp(t\nu) [A] ), \xi_C) = (\xi_C,\xi_C) + (e(t),\xi_C) .$$
Since $f$ is $K(S)$-invariant we also have 
$$ \Vert \xi_C \Vert_0^2  = \Vert \xi_C + t[\nu,\xi_C] + e(t) \Vert^2_0 $$
which implies 
$$ 2 (\xi_C,e(t)) = -t^2 \Vert [\nu,\xi_C] \Vert^2_0 + O(t^3) 
\ \text{as} \ t \to 0 .$$
Therefore,
$$ (\Phi(\exp(t\nu)[A]),\xi_C) = \Vert \xi_C \Vert^2_0 
- \hh t^2 \Vert [\nu,\xi_C] \Vert^2_0 + O(t^3) 
\ \text{as} \ t \to 0.$$
Since $T_{[A]}Y_C$ is the sum of the non-negative eigenspaces of the
Hessian of $(\Phi,\xi_C)$, it follows that $[\nu,\xi_C] = 0$.
\end{proof}

\blem \label{subs} There exists an open neighborhood $U$ 
of $C \cap Y_C^{\ss}$ in $Y_C^{\ss}$ such that if $ku \in U$ for some
$k \in K(S)$ and $u \in U$ then $k \in K_C(S)$.
\elem

\begin{proof} 
By the inverse function theorem, $K(S) \times_{K_C(S)} Y_C \to \M(X)$
is a local diffeomorphism onto its image in a neighborhood of $C$.
Hence there exists a neighborhood $U$ of $[A]$ in $Y_C$ and $V$ of
$K_C(S)$ in $K(S)$ such that for all $[A'] \in U$ and $k \in K(S)$, $
k \cdot [A'] \in Y_C$ implies $k \notin V$.  On the other hand, for
$\delta$ and $U$ sufficiently small, the set of $k \in K(S)$ such that
for all $[A'] \in U$,
$$( \Phi(k[A']), \xi ) \ge ( \Phi([A]),\xi ) - \delta $$
is contained in $V$.  If $k[A'] \in Y_C$ then $( \Phi(k[A']), \xi )
\ge ( \Phi([A]),\xi ) $.  This forces $k \in V$ and $k \in K_C(S)$.
\end{proof}

\blem $G(S) Y_C^{\ss}$ is a smooth embedded complex submanifold
diffeomorphic to $G(S) \times_{P_C(S)} Y_C^{\ss}$. \elem

\begin{proof}  Suppose that $ [A] \in Y_C^{\ss}$ and $g[A] \in
Y_C^{\ss}$ for some $g \in G(S)$.  By Lemma \ref{set} for any
neighborhood $U$ of $\Phinv(\xi_C)$ in $Y_C^{\ss}$, there exists an
element $p \in P_C(S)$ such that $[A'] = p[A]$ lies in $U$.  Since
$G(S) = K(S)P_C(S)$, $gp^{-1} = p'k$ for some $p' \in P_C(S)$ and $k
\in K(S)$.  Since $Y_C^{\ss}$ is $P_C(S)$-invariant, we have $k [A']
\in Y_C^{\ss} .$ By Lemma \ref{subs}, we can choose $U$ so that
$K_C(S)$ is a component of $\{ k \in K(S), k[A'] \in U \}$.  For $U$
sufficiently small, we obtain $k \in K_C(S)$, and $g = p'kp \in
P_C(S)$.

This shows that $G(S) \times_{P_C(S)} Y_C^{\ss} \to G(S) Y_C^{\ss}$ is a
bijection.  To show it is a diffeomorphism, we must show that the
condition
$$ \{ \xi \in \g(S), \ \xi_{[A]} \in T_{[A]} Y_C^{\ss} \} = \p_C(S)$$
holds for all $[A] \in Y_C^{\ss}$. This is open, and by the Lemma
above, it holds in a neighborhood of $Z_C \cap \Phinv(C)$.  The
condition is also invariant under $P_C(S)$, hence by Lemma \ref{open}
it holds everywhere.
\end{proof}

We wish to show that $G(S) Y_C$ is a minimising submanifold for $f$ in
the sense of Kirwan.  

\begin{lemma} \label{ortho}
The $L^2$-orthogonal subspace $(T_{[A]} G(S)Y_C)^\perp$ to $T_{[A]}
G(S)Y_C$ is a complement to $T_{[A]} G(S)Y_C$.
\end{lemma}

\begin{proof}  
Let $[A] \in C \cap Y_C$ and $\g(S)[A] \subset T_{[A]}$ the span of
the generating vector fields for $\g(S)$.  The $L^2$-orthogonal
subspace to $\g(S)[A]$ is a finite-dimensional complement to
$\g(S)[A]$ in $T_{[A]}$, by the proof of Lemma \ref{split} (c).  The
Lemma follows since $T_{[A]} G(S)Y_C$ contains $\g(S)[A]$.
\end{proof}

\blem The Hessian of $f$ is negative definite on $(T_{[A]}
G(S)Y_C)^\perp$, for any $[A] \in C$.  \elem

\begin{proof} We argue as in Kirwan \cite[p. 55]{ki:co}. 
Let $[A_t]$ be a path with $[A_0] = [A]$ and $\ddt [A_t] \in (T_{[A]}
G(S)Y_C)^\perp$.  Since $\ddt [A_t]$ is perpendicular to the
generating vector fields for the action,
$$ \Phi([A_t]) = \Phi([A]) + e(t) $$
where $e(t) = O(t^2)$.  Hence
$$ 2f([A_t]) = 2f([A]) + 2 \int_S \Tr(*_S \Phi([A]) \wedge e(t)) +
O(t^3) .$$
On the other hand,
$$ (\Phi([A_t]), \xi_C) = 2f([A]) + 2 \int_S \Tr(*_S \Phi([A]) \wedge
e(t)) .$$
It follows that the Hessians of $f$ and $(\Phi([A_t],\xi_C)$ agree up
to a scalar on $(T_{[A]} G(S) Y_C)^\perp$.  But $ (T_{[A]} G(S)
Y_C)^\perp$ is contained in $(T_{[A]} Y_C)^\perp $ on which the
Hessian of $(\Phi,\xi_C)$ is negative, by definition of $Y_C$.
\end{proof}

Together with \ref{set}, this shows that $G(S)Y_C$ is a minimising
manifold for $f$, in the sense of Kirwan.  It remains to show:

\begin{theorem} \label{flow} $G(S) Y_C^{\ss}$ is equal to $\M(X)_C$.
\end{theorem}

\begin{proof} We follow Kirwan \cite[p.91]{ki:co}.  Let $[A] \in C$ as
before and consider the splitting in Lemma \ref{ortho}.  By the
implicit function theorem there exists
\begin{enumerate}
\item  a neighborhood $V_A$ of $0$ in $\ker(d_A \oplus d_A^*)$, 
\item  a neighborhood $U_A$ of $[A]$ in $\M(X)$, and 
\item  a diffeomorphism $\varphi_A: \ U_A \to V_A$
\end{enumerate}
such that $U_A \cap G(S)Y_C$ is the pre-image of $V_A \cap T_{[A]}
G(S)Y_C$.  Let
$$H = \diag(H_{+ },H_{-})$$
be the decomposition of the Hessian of $f$ at $[A]$ into
positive-semi-definite and negative definite components.  In the slice
$U_A$ the trajectories $a_t = (a_{+,t},a_{-,t})$ are solutions to
$$ \ddt a_{+,t} = -H_+a_t + F_+(a_t)$$
$$ \ddt a_{-,t} = -H_{-}a_t + F_{-}(a_t) $$
where $F_+,F_{-}$ have vanishing derivatives at $(a_+, a_{-})
= (0,0)$.  It follows that
$$ a_{+,t} = e^{-H_+t} a_{+,0} + (\delta a_+)(t,a_0)  $$
$$ a_{-,t} = e^{-H_{-}t} a_{-,0} + (\delta a_{-})(t,a_0) $$
where $\delta a_+, \delta a_{-}$ have vanishing first partial
derivatives at the origin $a_0 = 0$.  For any $\eps > 0$ we may reduce
the size of $U_A$ so that
$$ (1 + \eps)^{-1} \Vert a_{-} \Vert \leq d(A,G(S)Y_C^{\ss}) \leq
(1+ \eps) \Vert a_{- } \Vert $$
everywhere in $U_A$.

\blem \label{blem} There exists a number $b > 1$, depending only on $C$,
such that if $U_A$ is taken sufficiently small then for any $a \in
U_A$ we have $\Vert a_{-,1} \Vert \ge b \Vert a_{-,0} \Vert.$
\elem

\begin{proof}  Let $c$ denote the minimum 
eigenvalue of $e^{-H_{-}}$ on $C$.  Choose $\theta > 0$ so that $c -
\theta > 1$, and $b = c - \theta$.  Since the partial Jacobian of
$\delta a_{-}$ vanishes at the origin, by shrinking $U_A$ we may
assume that $ \Vert \partial_{a_{-}} \delta a_{-} (1,a_+, a_{-}) \Vert
< \theta $ for all $a \in U_A$.  It follows that $
\Vert \delta a_{-} (1,a_+,a_{-}) \Vert \leq \theta \Vert a_{-} \Vert
$.  Hence for any $a_0 = (a_+, a_{-})$ we have
\begin{eqnarray*}
 \Vert a_{-,1} \Vert &=& \Vert e^{H_{-}} a_{-,0} 
+ \delta a_{-}(1,a_{ \ge,0},a_{-,0}) \Vert \\
&\ge& c \Vert a_{-,0} \Vert - \theta \Vert a_{-,0} \Vert \\
&\ge& b \Vert a_{-,0} \Vert .\end{eqnarray*}  
\end{proof}

\blem
 $\M(X)_C = G(S)Y_C^{\ss}$ in a neighborhood of $C$.  
\elem

\begin{proof} By 
Lemma \ref{blem} there is a neighborhood $U_C$ of $C$ in $\M(X)$ such
that if $[A_t]$ is a trajectory of \eqref{Aflow} then
$$d([A_1],G(S)Y_C^{\ss}) \ge b(1
+ \eps)^{-1}d([A_0], G(S)Y_C^{\ss}) .$$  
Choose $\eps$ sufficiently small so that $b(1 + \eps)^{-1} > 1$.  By
Lemma \ref{close} there exists a neighborhood $V_C$ of $C$ in
$\M(X)$ such that if $[A] \in V_C \cap \M(X)_C$ then $[A_t]$ lies in
$U_C$ for all $t \in [0,\infty)$.  Then for any $n \ge 1$,
$$ d([A_n],G(S)Y_C^{\ss}) \ge (b(1 + \eps)^{-1})^n d([A_0],
G(S)Y_C^{\ss}) .$$
But we may assume without loss of generality that
$d([A],G(S)Y_C^{\ss})$ is bounded on $U_C$.  Hence
$d([A],G(S)Y_C^{\ss}) = 0$.  This shows that $\M(X)_C \subset
G(S)Y_C^{\ss}$ in a neighborhood of $C$.  The opposite inclusion
follows from Lemma \ref{set}.
\end{proof} 

We now complete the proof of Theorem \ref{flow}.  Suppose $[A] \in
\M(X)_C$.  Then $[A_\infty] \in C$ and so the trajectory $[A_t]$
intersects $V_C$.  Since $[A_t] \in G(S)[A]$ this implies that $[A]
\in G(S)Y_C^{\ss}$.  Hence $\M(X)_C \subset G(S) Y_C^{\ss}$. By
Proposition \ref{set}, the subsets $G(S)Y_C^{\ss}$ are disjoint.
Since $\M(X)$ is the union of the stable manifolds $\M(X)_C$, we must
have $\M(X)_C = G(S)Y_C^{\ss}$.
\end{proof}

\begin{remark} \label{mb} Suppose that $f$ is a Morse-Bott function, that is,
the Hessian of $f |_C$ is non-degenerate along the normal bundle to
$S$. The map given by the time-1 flow $A_0 \mapsto A_1$ is hyperbolic,
using the estimates on the operator $L$ above.  It is then a
consequence of the stable manifold theorem \cite[Theorem III.8]{sh:gl}
that the strata $\M(X)_C$ are smooth.
\end{remark}

\section{Applications}

\subsection{K\"ahler quantization commutes with reduction}

Recall that $\L(X) \to \M(X)$ is the Chern-Simons pre-quantum line
bundle.  Let
$$\iota: \M(X)_{\Phinv(0)} \to \M(X), \ \ \ 
\pi: \ \M(X)_{\Phinv(0)} \to \M(\ol{X})$$
denote the inclusion, resp. the projection of the zero level set into
$\M(X)$, resp. onto $\M(\ol{X})$.
\begin{theorem} \label{QR}  Suppose $K(S)$ acts on $\Phinv(0)$ with finite
stabilizers.\footnote{This happens only if $X$ has markings.}
Then there is an isomorphism of global sections
$$H^0(\M(X),\L(X))^{G(S)} \to H^0(\M(\ol{X}),\L(\ol{X})), \ \ s
\mapsto \ol{s}, \ \ \iota^*s = \pi^* \ol{s} .$$
\end{theorem}

\begin{proof} The assumption on the stabilizers implies that 
the stable and semistable loci in $\M(X)$ coincide.  Therefore, for
any $[A] \in \M(X)^{\s}$ there exists an element $g_\infty$ such that
$[A_\infty] = g_\infty[A]$.  By Theorem \ref{semistable} (g) the image
$[g_\infty] \in G(S)/K(S)$ of $g_\infty$ is unique.  An application of
the implicit function theorem for Banach spaces shows that $g_\infty$
depends holomorphically on $[A]$.  Let $\L(X)^{\s}$ denote the
restriction of $\L(X)$ to $\M(X)^{\s}$, and $\L(X)_0$ the restriction
to $\Phinv(0)$.  For any $s \in H^0(\M(\ol{X},L(\ol{X}))$ define
$\ol{s} \in H^0( \M(X),\L(X)^{\s})^{G(S)}$ by $\hat{s}([A]) = s(
\pi_0( g_\infty[A]))$.  Since $\M(X)^{\s} \to \M(\ol{X})$ is
holomorphic, $\ol{s}$ is a holomorphic section. Also, $\Vert \ol{s}
\Vert$ is bounded by $\Vert s \Vert$.  By the Riemann extension
theorem for complex Banach manifolds (\cite[II.1.15]{ra:so}; a
convenient reference is \cite[Appendix]{ha:ba}) $\ol{s}$ has a
holomorphic extension to $\M(X)$.  Conversely, given $\ol{s}$ the
restriction $\iota^*\ol{s}$ is $K(S)$-invariant and so descends to
$\M(X)$.
\end{proof}
We remark that Teleman \cite{te:qu} has shown vanishing results for
the higher cohomology of the bundle $\L(\ol{X})$, using the
stratification of $\A(\ol{X})$ into Harder-Narasimhan types.
 
\subsection{Analog of Kirwan surjectivity}
\label{ksurj}

Let $H^\bullet_{K(S)}(\cdot, \Q)$ denote $K(S)$-equivariant cohomology
with rational coefficients.  If $K(S)$ acts with finite stabilizers on
$\Phinv(0)$, then $H^{\bullet}_{K(S)}(\Phinv(0),\Q) \cong
H^{\bullet}(\M(\ol{X}),\Q)$.

\begin{theorem} \label{Ksurj}
The inclusion $\Phinv(0) \to \M(X)$ induces a surjection
$H^{\bullet}_{K(S)} (\M(X),\Q) \to H^{\bullet}_{K(S)}(\Phinv(0),\Q)$.
\end{theorem}

\begin{proof}  This follows from the criterion of
Atiyah and Bott \cite[13.4]{at:mo}: For any critical component, the
circle $U(1)_C$ acts non-trivially on the normal bundle to $\M(X)_C$
at $Z_C$.  It follows that the Euler class for the normal bundle to
each stratum is invertible.  This implies that the stratification is
equivariantly perfect.
\end{proof}

The case that $S$ bounds a disk Theorem \ref{Ksurj} is a special case
of a recent result of Bott, Tolman, and Weitsman \cite{btw:ta}.  In
this case one has $X = X_+ \cup X_-$, say $X_+$ is homeomorphic 
to a disk, and 
$$ H^{\bullet}_{K(S)}( \M(X),\Q) = H^{\bullet}_{K(S)} (\M(X_-) \times \Omega K,\Q) =
H^{\bullet}_K(\M(X_-),\Q) $$
so the $H^\bullet_K(\M(X_-),\Q)$ surjects onto
$H^{\bullet}(\M(\ol{X}),\Q)$.  Their proof uses the Morse theory of
the energy functional on a space homotopy equivalent to $\M(X_-)$.
Inductive formulas for the Poincar\'e polynomials of the moduli spaces
of parabolic bundles using these techniques have been given by Nitsure
\cite{ni:co} and Holla \cite{ho:po}.

\begin{example}
To compare this approach with that of Atiyah and Bott, we compute the
Poincar\'e polynomial of the simplest case $\M(\ol{X};SU(2);\mu)$,
where $\mu = \hh$ is the marking corresponding to holonomy $-1$ around
a puncture.  Here and from now on, we identify $\t \to \R$ so that
$\Alc = [0,\hh]$.  Let $S$ be a circle enclosing the puncture. Then $X
= X_+ \cup X_-$, where $X_+$ is a surface of genus $2g$ with a single
boundary and no markings, and $X_-$ is a punctured disk.  By the
holonomy description Lemma \ref{surf} (\ref{holonomy})
$$ \M(X_+) = K^{2g} \times_K \Omega^1(S,\k),
\ \ \ \M(X_-) = K(S)/K(S)_\mu$$
where $K(S)_\mu$ denotes the stabilizer of $\mu$.  (Note this is
isomorphic to $K$ but does {\em not} consist of constant loops.)  Note
that $\M(X_+)$ is a principal $\Omega K$-bundle over $K^{2g}$.
Because the commutator $K^{2g} \to K$ induces the trivial map in
cohomology $H^{\bullet}(K) \to H^{\bullet}(K^{2g})$, the spectral
sequence for this fibration collapses at the second term.  Hence
$$ P_{K(S)}(\M(X)) = P_{K(S)_\mu} (\M(X_+)) = \frac{ P(\Omega K)
P(K^{2g})}{ 1- t^4} = \frac{(1 + t^3)^{2g}}{(1 - t^2)(1 - t^4)} .$$

Now we analyze the critical components.  Each critical component
contains a point with $\Phi(m_1,m_2) = *_S (\lambda,\mu)$, for some
$\lambda \in \Z_{> 0}$.  The corresponding critical component
$C_\lambda$ is the $K(S)$-orbit of $(T^{2g} \times \{ \lambda \})
\times \{\mu \}$ in the holonomy description of $ \M(X_+) \times
\M(X_-).$
The equivariant Poincar\'e polynomial of $C_\lambda$ is therefore
$$ P_{K(S)}(C_\lambda) = P_{U(1)}(T^{2g}) = \frac{(1+t)^{2g}}{1 - t^2}
.$$

\begin{lemma}  The index of $f$ at $C_\lambda$ is $g + 2\lambda - 2$.
\end{lemma}

\begin{proof}  The tangent space at
any point in $C_\lambda$ is isomorphic to 
$$ (\k/\t)^{2g} \oplus \k(S)/\k(S)_\lambda \oplus \k(S)/\k(S)_\mu.$$
The moment map to second order is
\begin{eqnarray*}
\Phi(t(\xi_1,\ldots,\xi_{2g},\zeta_-,\zeta_+))
&=&  \lambda + t^2 \sum_{i=1}^g [\xi_i,\xi_{g+i}] + t [\zeta_+,\lambda] + t^2
[\zeta_+,[\zeta_+,\lambda]]/2  \\
&-& \mu - t[\zeta_-,\mu] - t^2[\zeta_-,[\zeta_-,\mu]]/2 + O(t^3).\end{eqnarray*}
The second order term in $f$ is
\begin{multline*} 
 \int_S \Tr(\lambda - \mu, \sum_{i=1}^g [\xi_i,\xi_{i+g}]) 
+\frac{\Tr([\zeta_-,\lambda],[\zeta_-,\mu])+\Tr([\zeta_+,\lambda],[\zeta_+,\mu])}{2}
\\-\Tr([\zeta_+,\lambda],[\zeta_-,\mu]) . 
\end{multline*}
It follows that the index of $f$ on $(\k/\t)^{2g}$ is $g$.  Note that
if $\xi= 0 $ and $\zeta_- = \zeta_+$ then then the Hessian vanishes;
this are the directions tangent to the $K(S)$-orbit.  On the other
hand, suppose that $\zeta_- = - \zeta_+$ and is in the root space for
an affine root $\alpha$.  Then the Hessian is negative if and only if
$\alpha$ is negative on $\lambda$ and positive on $\mu$, or
vice-versa.  Therefore, the index of $f$ on $\k(S)/\k(S)_\lambda
\oplus \k(S)/\k(S)_\mu$ is the number of affine hyperplanes separating
$\lambda$ and $\mu$, equal to $2\lambda - 2$ for $\lambda > 0$.
\end{proof}

Putting everything together, the Poincar\'e polynomial of
$\M(\ol{X},\mu)$ is
\begin{eqnarray*}
 P(\M(\ol{X},\mu)) &=& P_{K(S)}(\M(X)) - \sum_{\lambda > 0}
t^{2(2\lambda + g - 2)} P_{K(S)}(C_\lambda) \\
&=&  \frac{(1 + t^3)^{2g} - (1+t)^{2g}t^{2g}}{(1 - t^2)(1 - t^4)}  
\end{eqnarray*}
which agrees with \cite[p. 335]{at:mo}.  
\end{example}

\subsection{Another proof of Birkhoff factorization} 

Let $\ol{X} = \P^1$ so that $\M(X) = \Omega K \times \Omega K$.  Since
the genus of $\ol{X}$ is zero this case is not covered by Theorem
\ref{main}.  However, $f$ is Morse-Bott which implies that the stable
manifolds are smooth, see Remark \ref{mb}.  The critical components
for $f$ are the orbits under $K(S)$ of pairs $(w_1\cdot 0,w_2 \cdot
0)$, and therefore can be indexed by $\Waff \backslash (\Waff \times
\Waff)/(W \times W)$.  The orbit $G(S) (w_1\cdot 0,w_2 \cdot 0)$ is a
submanifold with the same codimension as $\M(X)_{ [w_1,w_2]}$.  By
convergence of trajectories to critical points, any $G(S)$-orbit in
$\M(X)_{[w_1,w_2]}$ contains $K(S)(w_1\cdot 0,w_2 \cdot 0)$ in its
closure.  It follows that $\M(X)_{ [w_1,w_2]}$ is actually equal to
$G(S) (w_1\cdot 0,w_2 \cdot 0)$ and so
$$ \M(X) = \Omega K \times \Omega K = \bigcup_{ \Waff \backslash
(\Waff \times \Waff) / (W \times W)} G(S) (w_1\cdot 0,w_2 \cdot 0) .$$
Using $\M(X_-) = \Omega K = G(S)/G_{\hol}(X_-)$ we obtain (see
 \cite[Ch. 8]{ps:lg})
\begin{theorem}  $G(S) = G_{\hol}(X_-) \Lambda_+
G_{\hol}(X_+)$.\end{theorem}

\appendix

\section{Sobolev spaces}

A convenient reference for the basic material on Sobolev spaces is
\cite[Appendix]{ho:an3}.  The remaining results are variants of
results in Lions-Magenes \cite{lm:nh} which we learned from
R{\aa}de \cite[Section D]{rade:th}.

For any $s \in \R$, the Sobolev space $H^s(\R^n;\R^r)$ is the
completion of the space of smooth, compactly supported maps $u \in
C^\infty_0(\R^n;\R^r)$ in the norm
$$ \Vert u\Vert_s = (2 \pi)^{-n/2} \Vert (1 + \Vert \xi \Vert^2)^{s/2}
\mathcal{F}(u)(\xi) \Vert_{L^2(\R^n;\R^r)} $$
where $\mathcal{F}$ denotes Fourier transform.  Let $M$ be a compact
smooth manifold, and $E \to M$ a smooth Euclidean vector bundle of
rank $r$.  To define the Sobolev spaces of sections of $E$, we
introduce a trivializing atlas $\{ (U_i,V_i,\varphi_i: E |_{U_i} \to
V_i \times \R^r)\}$ and a subordinate partition of unity $\rho_i$, for
$i \in I$.  For any real $s$, define $H^s(M;E)$ to be the completion
of the space of smooth sections of $E$ with respect to the norm
$$ \Vert u \Vert_s = \left( \sum_{i \in I} \Vert \varphi_i \circ
(\rho_i u) \Vert_s^2 \right)^{\hh} .$$

Suppose that $M$ is a manifold with boundary, $\ti{M}$ is a closed
manifold of the same dimension containing $M$ (for instance, the
double of $M$), and $\ti{E} \to \ti{M}$ is a vector bundle whose
restriction to $M$ is isomorphic to $E$.  Let $r: \
C_0^\infty(\ti{M};\ti{E}) \to C^\infty_0(M;E)$ denote the restriction
map, and $\ol{C}_0^\infty(M;E)$ the image of $r$; this is the space of
{\em extendable} sections of $E$ with compact support.  The Sobolev
space $H^s(M;E)$ is the completion of $\ol{C}_0^\infty(M;E)$ in the
norm
$$ \Vert u \Vert_s = \inf_{\ti{u}} \Vert \ti{u} \Vert_s $$
where the infimum is over $\ti{u}$ that restrict to $u$.  Similarly,
$H^s_\partial(M;E)$ is the completion of the space $C^\infty_0(M
\backslash \partial M;E)$, with respect to the Sobolev norm induced by
the extension-by-zero map $C^\infty_0(M \backslash \partial M;E) \to
C^\infty_0(\ti{M})$.  These spaces have the following properties:
\blem \label{sobolev}
\ben
\item For $s > n/2$, the spaces $H^s(M;E) $ resp. $H^s_\partial(M;E)$
embed into the space $C(M;E)$ of continuous sections of $E$, resp.
vanishing on the boundary.
\item For any real $s$, there is a perfect pairing $ H^s(M;E) \times
H^{-s}_\partial(M;E) \to \R .$
\item When $s > \hh$, the elements of $H^s(M;E)$ have boundary values
(traces) in $H^{s- \hh}(\partial M;\partial E)$.  More generally,
choose a smooth vector field $\nu$ on $M$ normal to the boundary.  For
$m \in \N$ and $ s > m - \hh$, we have the Cauchy trace operator
$$ H^s(M) \to \prod_{j=0}^{m-1} H^{s - j - \hh}(\partial M), \ \ \ u
\mapsto \left( u |_{\partial M},\partial_{\nu} u|_{\partial M}, \ldots,
\partial_{\nu}^m u|_{\partial M} \right) $$
where $\partial_{\nu}$ denotes Lie derivative.
\item For $s \in (m- \hh, m + \hh)$, the kernel of the Cauchy trace
map is equal to $H^s_\partial(M)$.
\item \label{extension} 
For $s > \hh$, there is a continuous extension operator 
$$\mE: \ H^{s- \hh}(\partial M;\partial E) \to H^s(M;E)$$ 
which is a right inverse to $r_{\partial X}$, with the property that
for any $f \in H^{s- \hh}(\partial M;\partial E)$, $\mE(f)$ is smooth
away from $\partial M$.
\item For any real $s_1,s_2$, vector bundles $E_1,E_2 \to M$ and $s
\leq \min(s_1,s_2,s_1 + s_2 - n/2)$ (except for the borderline cases
$s_2 = -s_1$ and $s = -n/2$; $s = s_1$ and $s_2 = n/2$; $s = s_2$ and
$s_1 = n/2$) there is a continuous map
$ H^{s_1}(M;E_1) \times H^{s_2}(M;E_2) \to H^s(M;E_1 \otimes E_2) .$
\een
\elem

We will also need Sobolev spaces of mixed order on $[0,T] \times M$.
We assume throughout that $T \in (0,1)$.  Let $\ul{E}$ denote the
pullback of $E \to M$ to $[0,T] \times M$.  For any real $r,s$, the
space $ H^{r,s}([0,T] \times M;\ul{E})$ denotes the completion of the
space of smooth time-dependent sections of $E$ in the norm
$$ \Vert u \Vert_{r,s} = \inf_{\ti{u}} \Vert (\tau^2 + T^{-2})^{r/2}
\mathcal{F}(\hat{u})(\tau) \Vert_{L^2([0,T],H_s)} $$
where the infimum is over $\ti{u} \in C^\infty(\R,H^s(M;E))$ that
restrict to $u$ on $[0,T]$.  The space $ H^{r,s}_0([0,T] \times M;\ul{E}) =
H^{r,s}_0( [0,T] \times M;\ul{E}) $ is defined in the same way except that
the infimum is taken over $\ti{u} \in C^\infty(\R,H^s(M;E))$ that
restrict to $u$ on $[0,T]$ and vanish for $t < 0 $.  These spaces have
the following properties:
\blem \label{Sob}
\ben
\item For any real $r$ and $s$ the identity map defines an operator
$H_{0,r,s}([0,T] \times M; \ul{E}) \to H_{r,s}([0,T] \times
M;\ul{E})$.  For $r > -\hh$, this operator is injective. For $r < \hh$
it is onto.
\item For $r_1 \ge r_2$ and $s_1 \ge s_2$, the identity map defines an
operator $H^{r_1,s_1}([0,T] \times M;\ul{E}) \to H^{r_2,s_2}([0,T]
\times M;\ul{E})$ of norm less than $c T^{r_1 - r_2}$.
\item  For $r > \hh$ the map $f \mapsto f(0)$ extends to a 
restriction (trace) map $ H^{r,s}([0,T] \times M;\ul{E}) \to H^s(M;E)$.
More generally, for $m \in \N$, the Cauchy trace operator
\begin{multline}
 H^{r,s}([0,T] \times M;\ul{E}) \to \prod_{j=0}^{m-1} H^{r - j -
\hh,s}(M;E), \\ u \mapsto \left( u(0),u'(0),\ldots,u^{(m)}(0) \right) \end{multline}
is continuous for $ r < m - \hh$.  
\item For $r \in (m- \hh, m + \hh)$, the kernel of the Cauchy trace
map is equal to $H^{r,s}_0(M;E)$.
\item Let $E_1 \to M$ and $E_2 \to M$ be restrictions of vector
bundles $\ti{E}_1 \to \ti{M}, \ \ti{E}_2 \to \ti{M}$ to $M$.  For real
$r_1,s_1,r_2,s_2,r,s$ with $r \le \min(r_1,r_2,r_1+r_2 - \hh)$ and $s
\le \min(s_1,s_2,s_1 + s_2 - n/2)$ (except for the borderline cases
$r_2 = -r_1$ and $r = -1/2$, $r = r_1$ and $r_2 = 1/2$; $r = r_2$ and
$r_1 = 1/2$; $s_2 = -s_1$ and $s = -n/2$; $s = s_1$ and $s_1 = n/2$;
$s = s_2$ and $s_1 = n/2$.) there is a continuous map
$ H^{r_1,s_1}([0,T] \times M;\ul{E}_1) \times H^{r_2,s_2}([0,T]
\times M;\ul{E}_2) \to H^{r,s}([0,T] \times M;\ul{E}_1 \otimes
\ul{E}_2) .$
\item For any real $r,s$, integration with respect to Lebesgue measure
on $[0,T]$ defines a continuous map $ H_0^{r,s}([0,T] \times M;\ul{E})
\to H_0^{r+1,s}([0,T] \times M;\ul{E})$.  (Integration is not defined
on the spaces $H^{r,s}([0,T] \times M;\ul{E})$ for $r < - \hh$.)
\item For real numbers $r,s,r',s'$ the intersection $H^{r,s}([0,T]
\times M;E) \cap H^{r',s'}([0,T] \times M;\ul{E})$ is a Banach space with
norm $ \Vert \ \ \ \ \Vert_{(r,s) \cap (r',s')} = \left( \Vert \ \
\Vert_{(r,s)} + \Vert \ \ \Vert_{ (r',s')} \right)^{1/2} .$
Interpolation: For any $\theta \in [0,1]$, there is an embedding
$H^{r,s}([0,T] \times M;\ul{E}) \cap H^{r',s'}([0,T] \times M;\ul{E}) \to
H^{\theta r + (1- \theta)r', \theta s + (1-\theta)s'}([0,T] \times
M;E)$.
\item Let $P \in \Psi DO(M;E,E)$ be an self-adjoint, non-negative
elliptic pseudodifferential operator on $E$ of order $m$.  For any
real $r,s$, the operator $ \ddt + P $ defines an invertible operator $
H_{0}^{r+1,s}([0,T] \cap M;E) \cap H_0^{r,s+m}([0,T] \cap M;E) \to
H^{r,s}(M;E) .$
\item (Same assumptions on $P$) For any real $s,r$, solving the
homogeneous initial value problem $(\ddt + P)u = 0, u(0) = v$ defines
an operator $H^s(M;E) \to H^{\hh + r,s-mr}([0,T] \times M;\ul{E}), \ \
v \mapsto u$ of order bounded by $c \max(1,T^r)$.
\een
\elem

For short, we denote by $ \Omega^k(M;E)_s$,
resp. $\Omega^k(M;E)_{r,s}$, resp. $\Omega^k(M;E)_{0,r,s}$ the spaces
$H^s(M;\Lambda^k(T^*M) \otimes E)$, resp. $H^{r,s}([0,T] \times
M;\Lambda^k(T^*M) \otimes \ul{E}).$, resp. $H^{r,s}_0([0,T] \times
M;\Lambda^k(T^*M) \otimes \ul{E})$.

\section{Convergence in finite dimensions}

In this section we show that the gradient flow of minus the norm
square of the moment map converges for finite dimensional Hamiltonian
$K$-manifolds with proper moment maps.  This was previously obtained
by Duistermaat (unpublished).  

Let $K$ be a compact Lie group with Lie algebra $\k$.  Let $(\ , \ ):
\k \times \k \to \R$ denote an invariant inner product.  Let $M$ be a
Hamiltonian $G$-manifold with proper moment map $\Phi$ and invariant
compatible almost complex structure $J$, and
$$f = \hh(\Phi,\Phi) .$$
For any $C \subset \crit(f)$, let $U(1)_C$ denote the one-parameter
subgroup generated by $\Phi(C) \cap \t_+^*$, $M_C$ the stable manifold
for $C$, and $Z_C$ the component of $M^{U(1)_C}$ 
containing $C$.

\begin{theorem} \label{convergef} Any
trajectory $m(t)$ of $-\grad(f)$ has a unique limit $m(\infty)$ as $t
\to \infty$.  There exist constants $c,k$ and a time $T$ such that if
either (1) $m(\infty)$ is contained in the principal orbit-type
stratum for $Z_C$ or (2) $G_m = G_{m(\infty)}$ then
$d(m(t),m(\infty)) \leq c e^{-kt}$
for $t > T$; otherwise
$d(m(t),m(\infty)) \leq ct^{-1/2} .$
The map $M_C \to C, \ m \mapsto m(\infty)$ is a deformation retract.
\end{theorem}

The basic equality the controls the speed of the trajectory is
\begin{equation} \label{control} \frac{d}{dt} f(m(t)) = (- \grad(f) f)(m(t)) =
 - \Vert \grad(f) \Vert^2.\end{equation}
Since the change in $f(m(t))$ is finite, the integral of $\Vert
\grad(f) \Vert^2$ is finite.  The goal to remove the power $2$, in
order to conclude that the length of the trajectory is finite.

\begin{theorem}  \label{converge2}  Let $m'$ be a critical point for
$\phi$.  There exist constants $\eps_1,c_1,k_1 > 0 $ such that for any
$m \in M$ such that if $d(m,m') < \eps_1$, then either $\Phi(m(t)) <
\Phi(m')$ for some $t > 0 $ or there exists $m'' \in \crit(f)$ such
that $m(t) \to m''$ as $t \to \infty$ and $ f (m'') = f(m')$.
\end{theorem}

\begin{proof}  Assume  $\Phi(m(t)) \ge
\Phi(m')$ for all $ t \ge 0$.  For any $\gamma \in (0,1)$
\begin{equation} \label{gamma}
 \frac{d}{dt} (f(m(t)) - f(m))^{1 - \gamma} = (\gamma - 1) \Vert
\grad(f)(m(t)) \Vert^2  (f(m(t)) - f(m'))^{- \gamma} \ .
\end{equation} 
 by \eqref{control}.  Now we bound $\grad(f)$ from below in a
 neighborhood of $m'$.

\begin{lemma} \label{lb}
There are constants $\gamma \in (0,1)$, $c_2 > 0$ and $\eps_2 >0$ such that for any $m \in
M$ for which $ d(m,m') < \eps_2$ we have
$$ \Vert \grad(f)(m) \Vert \ge c_2 |  f(m) - f(m') |^\gamma.
$$
\end{lemma}

\begin{proof}  The proof relies on the local form theorem \cite{gu:no}.
Let $\g_{m'}^\circ \subset \g^*$ denote the annihilator of $\g_{m'}$.
Let $N$ denote the orthogonal complement of the tangent space to the
orbit $T_{m'} (K \cdot m')$.  Let $S$ denote the quotient of $N$ by
the kernel of the symplectic form restricted to $N$, called the {\em
symplectic slice} at $m'$.  Let $\omega_S$ denote the two-form on $S$,
and $\Phi_S$ the moment map.  Note that $\Phi_S$ is quadratic in $s
\in S$.  The $G$-manifold
$$ M' = G \times_{G_{m'}} (S \oplus (\g_{\Phi(m')}^*
\cap \g_{m'}^\circ)) .$$
has a Hamiltonian $G$-structure with moment map 
\begin{equation} \label{momentf} \Phi': \ [g,s,\nu] \mapsto g \cdot ( \phi_S(s) + \nu +
\Phi(m')) .\end{equation}
There exists an open subset $U$ of $m'$ in $M$ and a symplectomorphism
of $U$ with a neighborhood $U'$ of $[1,0,0]$ in $M'$.  One can choose
the symplectomorphism so that metric at $[1,0,0]$ is the direct sum of
a metric $g_S$ on $S$, and metric on $\g/ \g_{m'}$ and
$\g_{\Phi(m')}^* \cap \g_{m'}^\circ$ induced by the metric on $\g$.
For $\lambda \in \g_{\Phi(m)} \cap \g_m^\perp$,
$$ J_{m'} ( \lambda,s,0) = (0,J_Ss, \lambda) $$
where $J_S$ is the almost complex structure on $S$.  
 
Let $m \in U$ and let $[ g,s,\lambda] $ be its image in $U'$.  Let
$\xi = \Phi(m')$.  By \eqref{momentf}
\begin{eqnarray}
 f(m) - f(m') &=&  \hh \Vert \Phi_S(s) + \xi \Vert^2 + 
\hh \Vert \lambda \Vert^2 - \hh \Vert \xi \Vert^2  \\
&=& \hh \Vert \Phi_S(s) \Vert^2 + ( \Phi_S(s),\xi) + 
\hh \Vert \lambda \Vert^2 . \label{difff}
\end{eqnarray}
Using the local model we can trivialize the tangent bundle near $m'$
and approximate
$$ J_m = J_{m'} + e(m')$$
where $e(m) \to 0$ as $m \to m'$.
\begin{equation} \label{gradf} -\grad(f)(m) =
-J_m (\lambda,(\xi + \Phi_S(s)) \cdot s,0) = (0,-J_S ( \xi +
\Phi_S(s)) \cdot s, -\lambda) + e'(m)
\end{equation}
where $e'(m)$ denotes terms generated by the difference $J_m - J_{m'}$.
Let $S_0 \subset S$ denote the fixed point set of $G_m$, and $S_0
\oplus S_1$ the fixed point set of $\xi$, so that
$$ S = S_0 \oplus S_1 \oplus S_2 .$$
By \eqref{gradf}, $\Vert \grad(f) \Vert^2$ vanishes on $S_0$, vanishes
below degree $6$ on $S_1$, degree $2$ on $S_2$ and $\g_{\Phi(m)}^*
\cap \g_m^\circ$ and is positive away from $S_0$.  On the other hand,
by \eqref{difff} $f(m) - f(m')$ vanishes on $S_0$, vanishes below
degree $4$ on $S_1$ and degree $2$ on $S_2$ and $\g_{\Phi(m)}^* \cap
\g_m^\circ$.  It follows that
$ \Vert \grad(f) \Vert^2 \ge c |f(m) - f(m')|^{3/2} $
if the projection of $s$ onto $S_1$ is non-trivial, and
$ \Vert \grad(f) \Vert^2 \ge c |f(m) - f(m')| $
otherwise.  The lemma follows, taking $\gamma = 3/4$ resp. $1/2$ in
the first resp. second case.
\end{proof}

By Lemma \ref{lb} and \eqref{gamma} if $d(m(t),m') < \eps_2$ then
\begin{equation} \label{this}
- \frac{d}{dt} (f(m(t)) - f(m))^{1 - \gamma} \geq c_4 \Vert \grad(f)
\Vert .\end{equation}
Suppose $ d(m(t),m') < \eps_2 $ for $t \in [t_0,t_1]$.  By integrating
with respect to $t$ we obtain
\begin{equation} \label{int} \int_{t_0}^{t_1} \Vert \grad(f)(m(t)) \Vert
\d t \leq c_5 | f(m(t_0)) - f(m') |^{1- \gamma}.
\end{equation}
Suppose that $d(m(t),m') \leq \eps_2$ for all $ t \ge 0$.  By the
limit of \eqref{int} as $t_1 \to \infty$, the length of the gradient
flow line $m(t), t \in [t_0, \infty]$ approaches $0$ as $t_0 \to
\infty$.  By the Cauchy criterion, $m(t)$ converges to a critical
point $m''$ as $t \to \infty$.

To show that $d(m(t),m) \leq \eps_2$ for all $ t \ge 0$, suppose that
$s$ is the smallest number greater than $0$ such that $d(m(s),m') \ge
\eps_2$.  By smooth dependence of the flow on initial conditions,
there exists a constant $c$ such that $d(m(t),m') \leq c_7 d(m(0),m')
$ for $t \in [0,1]$.  We choose $\eps_1$ so that $c_7 \eps_1 <
\eps_2$, and assume $d(m(0),m') < \eps_1$.  Then $d(m(t),m') \leq
\eps_2, \ t \in [0,1]$ which implies $s > 1$.  By smooth dependence on
initial conditions, for any $t$ there is an $\eps_3 > 0$ and a
constant $c_3 > 0$ such that
\begin{equation} \label{that} d(m,m') < \eps_3 \implies
 f(m(t)) - f(m') \leq c_3 d(m,m')^2 .\end{equation}
If $d(m(t_0),m') < \min(\eps_3,\eps_1)$ then by \eqref{that}
and \eqref{int}
\begin{equation} \label{int2} \int_{t_0}^{t_1} \Vert \grad(f)(m(t))
\Vert \d t < c_6 d(m(t_0),m')^{2 (1 - \gamma)} .\end{equation}
The length of the gradient flow line between $1$ and $s$ is at least
$$ d(m(s),m(1)) \ge d(m(s),m') - d(m(1),m') \ge \eps_2 - c_7 \eps_1
.$$
By \eqref{int2} the length of the gradient flow line between $0$ and
$s$ is less than $c_8 \eps_1^{2(1-\gamma)}$, which is a contradiction
for $\eps_1$ sufficiently small.  This completes
the proof of Theorem \ref{converge2}. \end{proof}

Finally we compute the rate of convergence.
$$ \frac{d}{dt} | f(m(t)) - f(m') | = - \Vert \grad(f)(m(t)) \Vert^2
\leq - c | f (m(t)) - f(m') |^{2 \gamma} $$
implies 
$$ f(m(t)) - f(m') \leq ce^{-kt}$$ 
if $\gamma = 1/2$ or
$$ f(m(t)) - f(m') \leq ct^{-2}$$
if $\gamma = 3/4$.  Sufficient conditions for a point $m$ converging
exponentially to a point $m(\infty)$ in $C$ are (1) $G_m =
G_{m(\infty)} $ which implies that for some time $t$ such that $m(t)$
lies in the local model, the projection of $m(t)$ onto the slice $S$
is contained in $S_0 \oplus S_2$; and (2) $ S^\xi = S^{G_{m(\infty)}}
$ where $S$ is the symplectic slice at $m$, and $S^\xi,
S^{G_{m(\infty)}}$ are the fixed representations.  Since $S^\xi$ is
contained in the tangent space to $Z_C$, (2) holds if
$G_{m(\infty)}$ acts trivially on $T_{m(\infty)} Z_C$.
Equivalently, $m(\infty)$ is contained in the principal orbit-type
stratum of $Z_C$.  That $m \mapsto m(\infty)$ is a deformation
retract follows from smooth dependence on initial conditions, and is
left to the reader.  This completes the proof of Theorem
\ref{convergef}.

\def\cprime{$'$} \def\cprime{$'$} \def\cprime{$'$}
  \def\polhk#1{\setbox0=\hbox{#1}{\ooalign{\hidewidth
  \lower1.5ex\hbox{`}\hidewidth\crcr\unhbox0}}}



Revised February 24, 2004
\end{document}